\documentclass[10pt]{article}
\usepackage{fancyheadings, amsmath, graphicx, psfrag, color ,amsfonts}

\font\net=cmtt8
\font\fontemail=cmtt12
\font\fontauthors=cmcsc10 scaled \magstep1
\font\footsc=cmcsc10 at 8truept
 at 8truept

\newtheorem{Th}{Theorem}
\newtheorem{Def}[Th]{Definition}
\newtheorem{Cor}[Th]{Corollary}
\newtheorem{Lem}[Th]{Lemma}
\newtheorem{Prop}[Th]{Proposition}

\def\QED{\hfill\vrule height 1.5ex width 1.4ex depth -.1ex \vskip 10pt}

\newcommand{\masection}[1]{\section{#1}\setcounter{Th}{0}} 

\newtheorem{prerem}[Th]{Remark}
\newenvironment{rem}{\begin{prerem}\rm }{\end{prerem}}

\def\bigcapp#1#2{\displaystyle\bigcap _{\begin{array}{cc}\scriptstyle #1\\
\scriptstyle #2\end{array}}}

\def\D{\mathop{\rm D\!\!}\nolimits}

\def\Diag{\mathop{\rm Diag}\nolimits}
\def\dim{\mathop{\rm dim}\nolimits}
\def\End{\mathop{\rm End}\nolimits}

\def\Id{\mathop{\rm Id}\nolimits}

\def\ker{\mathop{\rm ker}\nolimits}

\def\summ#1#2{\sum _{\begin{array}{cc}\scriptstyle #1\\\scriptstyle #2\end{array}}}

\def\Sp{\mathop{\rm Sp}\nolimits}
\def\stirling2 #1#2{\left\{\begin{matrix} #1\\#2\end{matrix}\right\}}
\def\transp{\mathop{\ \!\!^t}\nolimits}
\def\uu{\mathop{\bf u}\nolimits}
\def\Vect{\mathop{\rm Span}\nolimits}

\font\dsrom=dsrom10 scaled 1200
\def \11{\textrm{\dsrom{1}}}

\def\g#1{\mathbb #1}

\begin{document}

\begin{center}
\LARGE{\bf An algebraic approach to P\'olya processes}
\end{center}

\vskip 10pt
\begin{center}
{\fontauthors
Nicolas Pouyanne
}
\end{center}

\begin{center}
D\'epartement de math\'ematiques\\
LAMA UMR 8100 CNRS\\ 
Universit\'e de Versailles - Saint-Quentin\\
45, avenue des Etats-Unis\\
78035 Versailles cedex\\
\end{center}

\begin{center}
{\fontemail
pouyanne@math.uvsq.fr}
\end{center}

\vskip 20pt
\begin{minipage}{11,5 truecm}
\begin{center}
{\bf Abstract}
\end{center}
P\'olya processes are natural generalization of P\'olya-Eggenberger urn models.
This article presents a new approach of their asymptotic behaviour {\it via}
moments, based on the spectral decomposition of a suitable finite difference
transition operator on polynomial functions.
Especially, it provides new results for \emph{large} processes
(a P\'olya process is called \emph{small} when $1$ is simple eigenvalue of
its replacement matrix and when any other eigenvalue has a real part
$\leq 1/2$;
otherwise, it is called large).

\end{minipage}

\vskip 10pt
\begin{minipage}{11,5 truecm}
\begin{center}
{\bf R\'esum\'e}
\end{center}
Les processus de P\'olya sont une g\'en\'eralisation naturelle des mod\`eles
d'urnes de P\'olya-Eggenberger.
Cet article pr\'esente une nouvelle approche de leur comportement
asymptotique {\it via} les moments, bas\'ee sur la d\'ecomposition spectrale
d'un op\'erateur aux diff\'erences finies sur des espaces de polyn\^omes.
En particulier, elle fournit de nouveaux r\'esultats sur les \emph{grands}
processus
(un processus de P\'olya est dit \emph{petit} lorsque $1$ est valeur propre
simple de sa matrice de remplacement et lorsque toutes les autres valeurs
propres ont une partie r\'eelle $\leq 1/2$ ;
sinon, on dit qu'il est grand).
\end{minipage}


\tableofcontents

\masection{Introduction}
\label{intro}

Take an urn (with infinite capacity) containing first finitely many
balls of $s$ different colours named $1,\dots ,s$.
This initial composition of the urn can be described by an $s$-dimensional
vector $U_1$, the $k$-th coordinate of $U_1$ being the number of balls of
colour $k$ at time $1$.
Proceed then to successive draws of one ball at random in the urn, any ball
being at any time equally likely drawn.
After each draw, inspect the colour of the ball, put it back into the urn
and add new balls following at any time the same rule.
This rule, summed up by the so-called \emph{replacement matrix}
$$
R=(r_{i,j})_{1\leq i,j\leq s}\in{\cal M}_s(\g Z)
$$
consists in adding (algebraically), for any $j\in\{1,\dots ,s\}$,
$r_{i,j}$ balls of colour $j$ when a ball of colour $i$ has been drawn.
In particular, a negative entry of $R$ corresponds to subtraction of
balls from the urn, when it is possible.
The \emph{urn process} is the sequence $(U_n)_{n\geq 1}$ of random vectors with
nonnegative integer coordinates, the $k$-th coordinate of $U_n$ being
the number of balls of colour $k$ at time $n$, {\it i.e.} after the
$(n-1)$-st draw.

Such urn models seem to appear for the first time in \cite{P-E}.
In 1930, in its original article
{\it Sur quelques points de la th\'eorie des probabilit\'es}
(\cite{Polya}),
G. P\'olya makes a complete
study of the two-colour urn process having a replacement matrix of the form
$S.\Id _2$, $S\in\g Z_{\geq 1}$.

We will only consider {\it balanced} urns.
This means that all rows of $R$ have a constant entries' sum, say $S$.
Under this assumption, the number of added balls is $S$ at any time, so
that the total number of balls at time $n$ is non random.
Furthermore, we will only consider replacement matrices having nonnegative
off-diagonal entries.
Any diagonal entries may be negative but subtraction of balls of a given
colour may become impossible.
In order to avoid this extinction, one classically adds an arithmetical
assumption to the column of any negative diagonal entry in $R$
(see Definition \ref{defiPolya} and related comments).
An urn process submitted to all these hypotheses will be called
{\bf P\'olya-Eggenberger}, in reference to the work of these authors.

\vskip 5pt
A P\'olya-Eggenberger urn process can be viewed as a Markovian
random walk in the first quadrant of $\g R^s$ with finitely many possible
increments (the rows of $R$), the conditional transition probabilities
between times $n$ and $n+1$ being linear functions of the coordinates of the
vector at time $n$.
This point of view leads to the following natural generalization:
we will name \emph{P\'olya process} such a random walk in $\g R^s$ with
normalized balance ($S=1$), even if it does not come from an urn process,
\emph{i.e.} even if $U_1$ and $R$ have non-integer values.
Note that a P\'olya process as it is defined just below looks very much like
a P\'olya-Eggenberger urn process, with the only difference that instead of
counting a number of balls, we deal with a positive real quantity $l_k(X_n)$
associated with each colour $k$ (corresponding to the ``number of balls'' of
this colour at time $n$), which gives the propensity to pick this colour at
the next step.
In this setting, $w_k$ is the vector in $\g R^s$ defined by the fact that,
when colour $k$ has been drawn, then for all $j\in\{ 1,\dots ,s\}$, one adds
$l_j(w_k)$ ``balls'' of colour $j$ to the urn.
P\'olya processes generalize P\'olya-Eggenberger urns only because
this propensity may be real-valued (see comments after
Definition~\ref{defiPolya}).

\begin{Def}
\label{defiPolya}
Let $V$ be a real vector space of finite dimension $s\geq 1$.
Let $X_1$,
$w_1,\dots ,w_s$ be vectors of $V$
and 
$(l_k)_{1\leq k\leq s}$ be a basis of linear forms on $V$ satisfying
the following assumptions:

\noindent
i- (initialization hypothesis)
\begin{equation}
\label{initialisation}
X_1\neq 0
{\rm ~~and~~}
\forall k\in \{1,\dots ,s\},~l_k(X_1)\geq 0
;
\end{equation}
\noindent
ii- (balance hypothesis)
for all $k\in\{1,\dots ,s\}$,
\begin{equation}
\label{hypobalance}
\sum _{j=1}^{s}l_j(w_k)=1;
\end{equation}
\noindent
iii- (sufficient conditions of tenability\footnote{
Some authors prefer the vocable \emph{viability} instead of tenability.
This last word has been chosen in reference to recent literature on the
subject.
})
for all $k,k'\in\{1,\dots ,s\}$,
\begin{equation}
\label{tenable}
\left\{
\begin{array}{lr}
\displaystyle k\neq k'\Longrightarrow l_k
(w_{k'})\geq 0,&(\ref{tenable}.a)\\
\displaystyle l_k(w_k)\geq0 {\rm ~~or~~}
l_k(X_1){\g Z}+\sum _{j=1}^{s}l_k(w_j){\g Z}=l_k(w_k){\g Z}.
&(\ref{tenable}.b)
\end{array}
\right.
\end{equation}
The (discrete and finite dimensional) {\bf P\'olya process} associated
with these data is the $V\!$-valued random walk $(X_n)_{n\in\g Z_{\geq 1}}$
with increments in the finite set $\{w_1,\dots ,w_s\}$,
defined by $X_1$ and the induction:
for every $n\geq 1$ and $k\in\{1,\dots ,s\}$,
\begin{equation}
\label{probaTransition}
{\rm Prob~}(X_{n+1}=X_n+w_k |X_n)=\frac {l_k(X_n)}{n+\tau _1-1}
\end{equation}
where $\tau _1$ is the positive real number defined by
\begin{equation}
\label{deftau1}
\tau _1=\sum _{k=1}^sl_k(X_1).
\end{equation}
\end{Def}

The process is defined on the space of all trajectories of
$X_1+\sum _{1\leq k\leq s}\g Z_{\geq 0}w_k$ endowed with the natural
filtration $({\cal F}_n)_{n\geq 0}$ where ${\cal F}_n$ is the $\sigma$-field
generated by $X_1,\dots ,X_n$.
It is Markovian\footnote{
The time-homogeneity of the process is more explicit when one reads condition
(\ref{probaTransition}) with denominator $\sum _kl_k(X_n)$
instead of $n+\tau _1-1$ (use Relation~(\ref{balance})).
}
and the transition conditional probabilities between
times $n$ and $n+1$ depend {\bf linearly} on the state at time $n$, as stated
in equations~(\ref{probaTransition}).
Conditions (\ref{initialisation}) and (\ref{hypobalance}) are necessary
and sufficient for the random
vector $X_2$ to be well defined by Relation~(\ref{probaTransition});
a readily induction shows the deterministic relation
\begin{equation}
\label{balance}
\forall n\geq 1,~~\sum _{k=1}^sl_k(X_n)=n+\tau _1-1.
\end{equation}
Condition (\ref{tenable}) suffices to
guarantee that the process is well defined, {\it i.e.} that the numbers
$l_k(X_n)$ do not become negative so that the process does not extinguish
as can be checked by an elementary induction.
The arithmetical assumption (\ref{tenable}.{\it b}), which has become classical
(compare with \cite{Gouet}, \cite{JansonFunctional}, \cite{FGP} for urns)
is equivalent to the following one:
$l_k(w_{k})$ is nonnegative, or it divides $l_k(X_1)$ and
all the $l_k(w_j)$ as real numbers.
Actually, if conditioned on non extinction, all the results about P\'olya
processes in this article remain valid when conditions (\ref{tenable})
are removed from the definition.

\vskip 5pt
P\'olya processes are natural generalizations of P\'olya-Eggenberger urns in
the following sense
(see \cite{BagchiPal}, \cite{Gouet}, \cite{FGP}, \cite{Puyhaubert}
for base references on P\'olya-Eggenberger urns).
Take a P\'olya-Eggenberger $s$-colour urn process having replacement
matrix $R$ and vector $U_1$ as initial composition;
let $S$ be the common sum of $R$'s rows, assumed to be nonzero.
The data consisting in taking the rows of $\frac 1SR$ as vectors $w_k$'s,
the coordinate forms as forms $l_k$'s and $X_1=\frac 1SU_1$ as initial vector
define a P\'olya process $(X_n)_n$ on $\g R^s$,
the random vector $X_n$ being $1/S$ times the $1\times s$ matrix $U_n$
whose entries  are the numbers of balls of different colours after $n-1$ draws.
We will name this process {\it standardized urn process}.
Conversely, if one considers the forms $l_k$ of a P\'olya process as being the
coordinate forms of $V$ (choice of a basis of $V$), the matrix
whose rows are the coordinates of the $w_k$'s satisfies all hypotheses
of a P\'olya-Eggenberger urn's replacement matrix with balance $S=1$, except
that its entries are not integers but real numbers.
This matrix will still be called \emph{replacement matrix} of the process.
Note that the balance property is expressed in Relation (\ref{hypobalance}).
The definition of P\'olya processes is readily stable after linear change of
coordinates, when urn processes do not have this property.

\vskip 5pt
The present text deals with P\'olya processes, so that all its results are
valid for P\'olya-Eggenberger urn processes.
Such a process being given, different natural questions arise:
what is the distribution of the vector at any time $n$?
Can the random vector be renormalized to get convergence?
What kind (and speed) of convergence is obtained?
What is the asymptotic distribution of the process?

\vskip 5pt
Since the work of P\'olya and Eggenberger, many authors have considered such
models, sometimes with more general hypotheses,
often with restrictive assumptions.
Direct combinatoric attacks in some particular cases were first intended
(\cite{Polya}, \cite{P-E}, \cite{Friedman} for example).
In the last years, they have been considerably refined by analytic
considerations on generating functions in low dimensions by much more general
methods (\cite{FGP}, \cite{Puyhaubert}).
A second approach was first introduced in \cite{AthreyaKarlin}
and developed in \cite{JansonFunctional} and \cite{JansonTriangular},
viewing such urns as multitype branching processes.
It consists in embedding the process in continuous time, using martingale
arguments and coming back to discrete time.
This method provides convergence results.
One can find in \cite{FGP}, \cite{Puyhaubert} and \cite{JansonFunctional}
good surveys and references on the subject.

\vskip 5pt
A P\'olya process will be called \emph{small} when $1$ is simple eigenvalue
of the replacement matrix $R$ and when every other eigenvalue of $R$ has
a real part $\leq 1/2$.
Otherwise, it will be said \emph{large}.

Under some assumptions of irreducibility on $R$, it is well known that if
$(X_n)_n$ is a small P\'olya process, a normalization 
$(X_n-nv_1)/\sqrt{n\log ^\nu n}$
converges in law to a centered Gaussian vector, $v_1$ being a deterministic
vector and $\nu$ a nonnegative integer that depends only on the conjugacy class
of $R$ - see \cite{JansonFunctional} for a complete statement of that fact.
In the case of reducible small processes, convergence in law after
normalization
has been shown for several families of processes in low dimensions;
this concerns for instance urns with triangular replacement
matrix~(\cite{FGP},~\cite{JansonTriangular},~\cite{Puyhaubert},
Example {\bf 2-} in Subsection~\ref{ex-subsec}).
Found limit laws in these studies are most often non normal.

In the case of large P\'olya processes, a suitable normalization of the random
vector $X_n$ leads to an almost sure asymptotics,
as shown in Theorems~\ref{lpss} and~\ref{lpnss}, main results of the paper.
These results do not require any irreducibility assumption.
This asymptotics is described by finitely-many random variables $W_k$
that appear as limits of martingales.
Joint moments of the $W_k$ are computed in terms of so-called reduced
polynomials $(Q_\alpha )_{\alpha\in (\g Z_{\geq 0})^s}$ that will be defined
later and initial conditions of the process.
We give hereunder a simplified version of the result:
suppose that the replacement matrix $R$ has $1$ and $\lambda _2$ as
\emph{simple} eigenvalues and that any other eigenvalue is the conjugate
$\overline{\lambda _2}$ or has a real part $<Re (\lambda _2)$.
Such a process will be called \emph{generic}
\footnote{
Note that such a process \emph{is} generic in the sense that almost all
(in the strong sense of algebraic geometry) replacement matrices of P\'olya
processes satisfy this assumption.}.

\vskip 5pt\noindent
{\bf Asymptotics of generic large P\'olya processes}
\it
If $(X_n)_n$ is a generic large P\'olya process,
there exist some complex-valued random variable $W$ and non random complex
vectors $v_1$ and $v_2$ such that
$$
X_n=nv_1+\Re \left( n^{\lambda _2}Wv_2\right)
+o\left( n^{\Re(\lambda _2)}\right) ,
$$
the small $o$ being almost sure and in any ${\rm L}^p$, $p\geq 1$.
Furthermore, any joint moment of the variables $W$ and its complex
conjugate $\overline W$ is given by the formula
$$
E\left( W^p\overline{W^q}\right) =
\frac{\Gamma (\tau _1)}{\Gamma (\tau _1+p\lambda _2+q\overline{\lambda _2})}
Q_{(0,p,q,0,\dots )}(X_1)
$$
where $\Gamma$ is Euler's function.
\rm

\vskip 5pt\noindent
The positive number $\tau _1$, defined by~(\ref{deftau1}), depends on initial
condition $X_1$.
Vectors $v_1$ and $v_2$ are here eigenvectors of the replacement matrix
respectively associated with the eigenvalues $1$ and $\lambda _2$.
In particular, the second order term is oscillating when $\lambda _2$ is
non real, giving a complete answer to the already observed non convergence
of any non trivial normalization $(X_n-EX_n)/n^z$,
$z\in\g C$~(see \cite{ChernHwang} and related papers for example).

\vskip 5pt
The method used here to establish the general asymptotics of large P\'olya
processes also leads to results on distributions at finite time
(exact expressions for moments for example) but we do not focus
on this point of view.
It relies on asymptotic estimates of suitable moments of $X_n$.
Hence, the first step is to express, for general functions $f$, the expectation
$Ef(X_n)$ in terms of initial condition $X_1$ and of iterations of a finite
difference operator $\Phi$, namely, by Proposition~\ref{esper},
$$
Ef(X_n)=\gamma _{\tau _1,n}(\Phi)(f)(X_1)
$$
where $\gamma _{\tau _1,n}$ is the polynomial defined by
$\gamma _{\tau _1,1}=1$ and, for any $n\geq 2$,
\begin{equation}
\label{gamma}
\gamma _{\tau _1,n}(t)
=\prod _{k=1}^{n-1}\left( 1+\frac t{k+\tau _1-1}\right) ;
\end{equation}
$\Phi$ is the {\bf transition operator} associated with the process,
defined on the space of all functions
$f:V\to \g R$
(or more generally on the space of all functions $f:V\to W$ where $W$ is any
real vector space) by:
$\forall v\in V$,
\begin{equation}
\label{phi}
\Phi (f)(v)=
\sum _{1\leq k\leq s}l_k(v)\bigg[ f(v+w_k)-f(v)\bigg].
\end{equation}
The second step is to study this linear operator $\Phi$ on its restriction
to the space of linear forms on $V$, which leads to set a corresponding
Jordan basis $(u_k)_{1\leq k\leq s}$ of this space, with corresponding
eigenvalues $(\lambda _k)_{1\leq k\leq s}$
(Definition~\ref{jordanBasis-def}).
The third step consists in observing, as done in Proposition~\ref{phiordre},
that $\Phi$ stabilizes, for any
$\alpha\in (\g Z_{\geq 0})^s$, the finite dimensional polynomial subspace
$S_\alpha=\Vect\{\uu ^\beta ,~\beta\leq\alpha\}$
where, for all $\beta =(\beta _1,\dots ,\beta _s)\in (\g Z_{\geq 0})^s$,
$\uu ^\beta =\prod _{1\leq k\leq s}u_k^{\beta _k}$
and $\leq$ is the degree-antialphabetical order on $s$-uples of integers,
defined below by~(\ref{ordre}).
Therefore, it is subsequently possible to decompose any $\uu$-monomial
$\uu^\alpha$, $\alpha\in (\g Z_{\geq 0})^s$ as a sum of functions in the
characteristic subspaces\footnote
{When the context is unambiguous, if $z$ is a complex number, $z$ will
also denote $zI$ where $I$ is the identity endomorphism.
}
$\ker (\Phi-z)^\infty =\bigcup _{n\geq 0}\ker (\Phi-z)^n$, $z\in\g C$.

If one denotes $\lambda =(\lambda _1,\dots ,\lambda_s)$ and
$\langle\alpha ,\lambda\rangle =\sum _{1\leq k\leq s}\alpha _k\lambda _k$
for any $\alpha\in (\g Z_{\geq 0})^s$,
it turns out that the eigenvalues of the restriction of $\Phi$ to stable finite
dimensional polynomial spaces are precisely the
$\langle\alpha ,\lambda\rangle$, as justified in Section~\ref{results-sec}.
The projection of any $\uu ^\alpha$ on
$\ker (\Phi -\langle\alpha ,\lambda\rangle )^\infty$
parallel to
$\bigoplus _{z\neq\langle\alpha ,\lambda\rangle}\ker (\Phi -z)^\infty$
will be denoted by $Q_\alpha$ and named \emph{reduced polynomial} of $\Phi$ of
rank $\alpha$.
The reduced polynomials of rank $\leq\alpha$ constitute a basis of $S_\alpha$
and any $\uu ^\alpha$ can be written
\begin{equation}
\label{relationsu/Q}
\uu ^\alpha =Q_\alpha
+\sum _{\beta <\alpha ,~\langle\beta ,\lambda\rangle
\neq\langle\alpha ,\lambda\rangle}q_{\alpha ,\beta}Q_\beta
\end{equation}
as proved in Proposition~\ref{proprietesQalpha}.

This leads to an asymptotic estimate of the moments $E\uu ^\alpha (X_n)$
(Theorem~\ref{jointmoments})
since, for any $z\in\g C$ and any $f\in\ker (\Phi-z)^\infty$, there exists
an integer $\nu\geq 0$ such that
\begin{equation}
\label{formEf(Xn)}
Ef(X_n)\mathop{\sim}_{n\to +\infty}
\frac{n^z\log ^\nu n}{\nu !}\frac{\Gamma (\tau _1)}{\Gamma (\tau _1+z)}
(\Phi -z)^\nu (f)(X_1)
\end{equation}
as it is proven in Corollary~\ref{fjordan}.
The asymptotic estimate in Theorem~\ref{jointmoments} is based on the
determination of the indices $\beta$ in expansion~(\ref{relationsu/Q})
that contribute to the leading term of $E\uu ^\alpha (X_n)$;
this is the object of the whole
Subsections~\ref{cones-subsec} and ~\ref{continued-subsec}.
To this end, Theorem~\ref{coefflemma} enables to refine
Relation~(\ref{relationsu/Q}):
it implies that a coefficient $q_{\alpha ,\beta}$ does not vanish only if
$\beta$ belongs to a convex polyhedron
$(A_\alpha -\Sigma )\cap (\g R_{\geq 0})^s$ of $\g R^s$,
where $A_\alpha$ is a the set of nonnegative integer points of a certain
rational cone with vertex $\alpha$ that depends
on the P\'olya process and $\Sigma$ a universal rational cone
(universal means here that $\Sigma$ is the same one for any P\'olya process).
Definitions of $\Sigma$ and $A_\alpha$ are respectively given
by~(\ref{defSigma}) and~(\ref{formulaAalpha}).
Formula~(\ref{relationsu/Q}) can thus be refined into
\begin{equation}
\label{refinedRelationsu/Q}
\uu ^\alpha =Q_\alpha
+\sum _{\beta\in A_\alpha -\Sigma ,~\langle\beta ,\lambda\rangle
\neq\langle\alpha ,\lambda\rangle}q_{\alpha ,\beta}Q_\beta
\end{equation}
which is the same as Relation~(\ref{refined}).

We will say that
$\alpha =(\alpha _1,\dots ,\alpha _s)\in (\g Z_{\geq 0})^s$ is a
\emph{power of large projections} whenever
$\alpha _k=0$ for all indices $k$ such that $\Re (\lambda _k)\leq 1/2$;
similarly, $\alpha$ will be called \emph{power of small projections} whenever
$\alpha _k=0$ for all indices $k$ such that $\Re (\lambda _k)>1/2$.
Now, if $\alpha$ is a power of large projections,
Propositions~\ref{majorationsRe} {\it 1-} and~\ref{propAalphalambda}
imply that
$\Re\langle\beta ,\lambda\rangle <\Re \langle\alpha ,\lambda\rangle$
whenever $\beta\in A_\alpha -\Sigma ,~\langle\beta ,\lambda\rangle
\neq\langle\alpha ,\lambda\rangle$.
Therefore, thanks to Relation~(\ref{formEf(Xn)}), the leading term
of $E\uu ^\alpha (X_n)$ in Formula~(\ref{refinedRelationsu/Q}) will come
from $EQ_\alpha (X_n)$ only, with an order of magnitude of the form
$n^{\langle\alpha ,\lambda\rangle}\log ^\nu n$,
the number $\Re\langle\alpha ,\lambda\rangle$ being $>|\alpha |/2$.
Similarly, Propositions~\ref{majorationsRe} {\it 2-} implies that,
if $\alpha$ is a power of small projections, this order of
magnitude never exceeds $n^{|\alpha |/2}\log ^\nu n$ for some nonnegative
integer~$\nu$.
A precise statement of these moments' asymptotics is given in
Theorem~\ref{jointmoments}.
Note that the intervention of $\Sigma$ can be bypassed by a self-sufficient
argument that have been suggested by the anonymous referee
(see Remark~\ref{GammaReferee}).

\vskip 10pt
Section~\ref{Preliminaries} is devoted to Jordan decomposition of $\Phi$'s
restriction to linear forms and related definitions and notations.
The main results of the paper are introduced and completely stated in
Section~\ref{results-sec} while
the action of transition operator $\Phi$ on polynomials is studied in
Section~\ref{phi-sec}.
This is done in three steps: 
first, the stability of the filtration $(S_\alpha )_\alpha$ of subspaces is
established as well as its consequences on reduced polynomials;
cone $\Sigma$ and polyhedra $A_\alpha$ are then
introduced in the space $(\g R_{\geq 0})^s$ of exponents;
afterwards, consequences of these geometrical considerations are drawn to
refine $\Phi$'s action.
Main Theorems~\ref{jointmoments}, \ref{lpss} and~\ref{lpnss} are proved in
Sections~\ref{jointmoments-sec} and~\ref{main-sec}.
At last, Section~\ref{Exemples-sec} contains diverse remarks and examples.

\vskip 10pt\noindent
{\bf Acknowledgements}

The author is very grateful to the anonymous referee whose careful reading
of a preprint version of this paper and helpful comments on it greatly improved
its quality.

\masection{Preliminaries, notations and definitions}
\label{Preliminaries}

Definition of P\'olya processes in a real vector space $V$ of finite dimension
$s\geq 1$ was given in Definition~\ref{defiPolya}.
We associate with any process its \emph{replacement endomorphism}
that will be denoted by $A$ in reference to literature on the subject
(\cite{AthreyaKarlin}, \cite{JansonFunctional} for example).
Let $V_{\g C}=V\otimes _{\g R}\g C$ be the complexified space of $V$.

\begin{Def}
If $(X_n)_n$ is a P\'olya process, its {\bf replacement endomorphism} is,
with notations of Definition~\ref{defiPolya}, the endomorphism
$A=\sum _{1\leq k\leq s}l_k\otimes w_k\in V^*\otimes V\simeq\End (V)$,
defined as
$$
A(v)=\sum _{1\leq k\leq s}l_k(v)w_k
$$
for every $v$ in $V$.
\end{Def}

\noindent
Note that the transpose of $A$ is the restriction of the transition operator
$\Phi$ to linear forms on $V$.
When the process is a P\'olya-Eggenberger urn process, the matrix of $A$ in
the dual basis of $(l_k)_k$ is the transpose of the normalized urn's
replacement matrix $\frac 1SR$ (notations of Section~\ref{intro}).

With this definition, the expectation of $X_{n+1}$ conditionally to $X_n$ is
readily expressed as $(I+\frac A{n+\tau _1-1})X_n$, so that the expectation
of $X_n$ equals
$$
EX_n=\gamma _{\tau _1,n}(A)(X_1)
$$
(straightforward induction).

One of the first tools used to describe the asymptotics of a P\'olya process
is the reduction of its replacement endomorphism $A$
(or of its transpose on the dual vector space of $V$).
Because of condition (\ref{hypobalance}),
the linear form $u_1=\sum _{k=1}^sl_k$ satisfies $u_1\circ A=u_1$, which
shows that {\bf $1$ is always eigenvalue of $A$}.
The whole assumptions (\ref{initialisation}),(\ref{hypobalance})
and (\ref{tenable}), allows us to say more on $A$'s spectral decomposition.
Even if these properties can be proved using Perron-Frobenius theory, we
give a proof's hint of Proposition \ref{spectre}.

\begin{Prop}
\label{spectre}
Any complex eigenvalue $\lambda$  of $A$ equals $1$ or satisfies
$\Re\lambda <1$.
Moreover, $\dim\ker (A-1)$ equals the multiplicity of $1$ as eigenvalue of $A$.
\end{Prop}

\begin{pff}
Replace $A$ by its matrix in the dual basis of $(l_k)_k$.
Suppose first that all entries of $A$ are nonnegative.
The space of all $s\times s$ matrices having nonnegative entries and columns
with entries' sum $1$ is bounded (for the norms' topology) and stable for
multiplication.
This forces the sequence $(A^n)_{n\geq 0}$ to be bounded, which implies both
results (for the second one, consider Jordan's decomposition of $A$ and
note that the positive powers of  $I+N$ constitute an unbounded sequence
if $N$ is a nilpotent nonzero matrix).
If $A$ has at least one negative diagonal entry, apply the results to
$(A+a)/(1+a)$ for any positive $a$ such that $A+a$ has nonnegative entries.
\end{pff}\QED

In the whole paper, a P\'olya process with replacement endomorphism $A$ being
given, we will denote by $\sigma _2$ the real number $\leq 1$ defined by
\begin{equation}
\label{defsigma2}
\sigma _2=\left\{
\begin{array}{l}
1{\rm ~if~}1{\rm ~is~multiple~eigenvalue~of~}A;\\~\\
\max\{\Re\lambda ,~\lambda\in\Sp (A),~\lambda\neq 1\}{\rm ~otherwise},
\end{array}
\right.
\end{equation}
where $\Sp (A)$ is the set of eigenvalues of $A$.

\subsection{Jordan basis of linear forms of the process}
\label{jordanbasis}

The present subsection is devoted to notations and vocabulary related to
spectral properties of the replacement endomorphism $A$.

\begin{Def}
\label{jordanBasis-def}
If $(X_n)_n$ is a P\'olya process of dimension $s$, a basis
$(u_k)_{1\leq k\leq s}$ of linear forms on $V_{\g C}$ is called a
{\bf Jordan basis of linear forms of the process} or shortlier
a {\bf Jordan basis}
when

\noindent 1-
$u_1=\sum _{1\leq k\leq s}l_k$;

\vskip 3pt
\noindent 2-
$u_k\circ A=\lambda _ku_k+\varepsilon _ku_{k-1}$
for all $k\geq 2$, where
the $\lambda _k$ are complex numbers (necessarily eigenvalues of $A$) and
where the $\varepsilon _k$ are numbers in $\{ 0,1\}$
that satisfy
$\lambda _k\neq\lambda _{k-1}\Longrightarrow \varepsilon _k=0$.
\end{Def}
In other words, the matrix of the transposed endomorphism $^t\! A$
in a Jordan basis of linear forms
has a block-diagonal form
$\Diag \left( 1,J_{p_1}(\lambda _{k_1}),\dots ,J_{p_t}(\lambda _{k_t})\right)$
where $J_p(z)$ denotes the $p$-dimensional square matrix
$$
J_p(z)=
\left(
\begin{array}{cccc}
z&1&&\\
&z&\ddots &\\
&&\ddots &1\\
&&&z\\
\end{array}
\right) .
$$
A (real or complex) linear form $u_k$ will be called {\bf eigenform}
of the process when $u_k\circ A=\lambda _ku_k$, {\it i.e.} when
$\varepsilon _k=0$.
An eigenform of the process is an eigenvector of $^t\! A$;
some authors call these linear forms left eigenvectors of $A$, refering
to matrix operations.

\begin{Def}
A Jordan basis of linear forms being chosen with notations as above,
a subset $J\subseteq\{ 1,\dots ,s\}$ is called a
{\bf monogenic block of indices}
when $J$ has the form $J=\{m,m+1,\dots ,m+r\}$
($r\geq 0$, $m\geq 1$, $m+r\leq s$)
with $\varepsilon _m=0$, $\varepsilon _k=1$ for every
$k\in\{m+1,\dots ,m+r\}$ and $J$ is maximal for this property.
Any monogenic block of indices $J$ is associated with a unique eigenvalue
of $A$ that will be denoted by $\lambda (J)$.
\end{Def}
In other words, $J$ is monogenic when the subspace $\Vect\{ u_j,~j\in J\}$ is
$A$-stable and when the matrix of the endomorphism of $\Vect\{ u_j,~j\in J\}$
induced by $^t\! A$  in the Jordan basis is one of the Jordan blocks
mentioned above with number $\lambda (J)$ on its diagonal.
The adjective {\it monogenic} has been chosen because this means that the
subspace $\Vect\{u_j,~j\in J\}=\g C[^t\! A].u_{m+r}$ is a monogenic
sub-$\g C[t]$-module of the dual space $V_\g C^*$ for the usual
$\g C[t]$-module structure induced by $^t\! A$.

\begin{Def}
\label{defPrincipalBlock}
A monogenic block of indices $J$ is called a {\bf principal block} when
$\Re\lambda (J)=\sigma _2$ and $J$ has maximal size among the monogenic blocks
$J'$ such that $\Re\lambda (J')=\sigma _2$
(see~(\ref{defsigma2}) for $\sigma _2$'s definition).
\end{Def}

\noindent
A Jordan basis $(u_k)_{1\leq k\leq s}$ of linear forms of the process
being chosen,
\begin{equation}
\label{vk}
(v_k)_{1\leq k\leq s}
\end{equation}
will denote its dual basis, made of the vectors of $V_{\g C}$ that satisfy
$u_k(v_l)=\delta _{k,l}$ (Kronecker notation) for any $k$ and $l$,
and
\begin{equation}
\label{lambda}
\lambda =(\lambda _1,\dots ,\lambda _s)
\end{equation}
the $s$-uple of eigenvalues (distinct or not) respectively associated with
$u_1,\dots ,u_s$ (or $v_1,\dots ,v_s$).
In particular, $\lambda _1=1$ for any Jordan basis of linear forms.
The eigenvalues $\lambda _1,\dots ,\lambda _s$ of $A$ are called
{\bf roots of the process}.
For any $k$, we also denote by $\pi _k$ the projection on the line $\g Cv_k$
relative to the decomposition
$V_{\g C}=\bigoplus _{1\leq l\leq s}\g Cv_l$;
these projections satisfy
\begin{equation}
\label{projections}
\Id =\sum _{1\leq k\leq s}\pi _k
{\rm ~~and~}
\pi _k=u_k.v_k.
\end{equation}
Note that the $\pi _k$ commute with each other
($\pi _k\pi _l=\delta _{k,l}\pi _k$) but do not commute with $A$.
Nevertheless, $A$ commutes with $\sum _{j\in J}\pi _j$, the sum being
extended to any monogenic block of indices $J$ (these sums are polynomials
in $A$).
This fact will be used in the proofs of Theorems \ref{lpss} and \ref{lpnss}.
The lines spanned by the vectors $v_k$ can be seen as principal directions
of the process, the word principal being here used in physicists' sense.

\subsection{Semisimplicity, large and small projections}
\label{powers}

For every Jordan basis $(u_k)_{1\leq k\leq s}$ of linear forms,
and for every $\alpha =(\alpha _k)_{1\leq k\leq s}\in\g Z^s$,
we adopt the notations
\begin{equation}
\label{pdtscal}
\begin{array}{c}
\displaystyle |\alpha |=\sum _{1\leq k\leq s}\alpha _k
\hskip 20pt{(total~degree)}\\
\\
\displaystyle
\langle\alpha ,\lambda \rangle=\sum _{1\leq k\leq s}\alpha _k\lambda _k
\end{array}
\end{equation}
and, when all the $\alpha _k$ are nonnegative integers
\begin{equation*}
{\uu}^\alpha =\prod _{1\leq k\leq s}u_k^{\alpha _k},
\end{equation*}
${\uu}^\alpha$ being a homogeneous polynomial function of degree $|\alpha |$.

\vskip 5pt
Given a Jordan basis $(u_k)_{1\leq k\leq s}$ of linear forms of the process,
we adopt the following definitions.

\begin{Def}
A P\'olya process is called {\bf semisimple} when its replacement
endomorphism $A$ is semisimple, {\it i.e.} when $A$ admits a basis of
eigenvectors in $V_{\g C}$ (this means that all the $u_k$ are real or
complex eigenforms of $A$).
The process is called {\bf principally semisimple} when all principal
blocks have size one (for any choice of a Jordan basis).
\end{Def}
The four following assertions are readily equivalent:

\noindent
i) the process is principally semisimple;

\noindent
ii) for any $k\in\{1,\dots ,s\}$,
$\big(\Re\lambda _k=\sigma _2\Longrightarrow u_k$ {is eigenform}$\big)$;

\noindent
iii) the induced endomorphism
$\big(\sum _{\{k,~\Re\lambda _k=\sigma _2\}}\pi _k\big)A$
is diagonalizable over $\g C$;

\noindent
iv) if $r\geq 1$ and if $\{\lambda _k,~k\geq r+1\}$ are the roots of
the process having a real part $<\sigma _2$,
the matrix of $^t\! A$ in the Jordan basis
has a block-diagonal form
$\Diag \left( 1,\lambda _2,\dots ,\lambda _r,
J_{p_1}(\lambda _{k_1}),\dots ,J_{p_t}(\lambda _{k_t})\right)$.

\noindent
Note that Proposition \ref{spectre} asserts that any $u_k$ associated with
root $1$ is eigenform of $A$.

\begin{Def}
A root of the process is called {\bf small} when its real part
is $\leq 1/2$;
otherwise, its is said {\bf large}.
The process is called {\bf small} when $\sigma _2\leq 1/2$, which means
that $1$ is simple root and all other roots are small;
when the process is not small, it is said {\bf large}.
\end{Def}

\begin{Def}
\label{largesmall}
Let $\alpha =(\alpha _1,\dots ,\alpha _s)\in (\g Z_{\geq 0})^s$.

\noindent
1- $\alpha$ is called {\bf power of large projections} when $\uu ^\alpha$ is
a product of linear forms associated with large roots, {\it i.e.} when for
all $k\in\{ 1,\dots ,s\}$,
$\big( \alpha _k\neq 0\Longrightarrow\Re\lambda _k>1/2\big)$.

\noindent
2- $\alpha$ is called {\bf power of small projections} when $\uu ^\alpha$ is
a product of linear forms associated with small roots, {\it i.e.} when for
all $k\in\{ 1,\dots ,s\}$,
$\big( \alpha _k\neq 0\Longrightarrow\Re\lambda _k\leq 1/2\big)$.

\noindent
3- $\alpha$ is called {\bf semisimple power} when $\uu ^\alpha$ is a
product of eigenforms, {\it i.e.} when for all $k\in\{ 1,\dots ,s\}$,
$\big( \alpha _k\neq 0\Longrightarrow u_k$
{\it is eigenform of the process}$\big)$.

\noindent
4- $\alpha$ is called {\bf monogenic power} when its support in contained in
a monogenic block of indices.
\end{Def}

\vskip 5pt
In the whole text, the canonical basis of $\g Z^s$ (or of $\g R^s$)
will be denoted by
\begin{equation}
\label{deltak}
(\delta _k)_{1\leq k\leq s}
\end{equation}
and the symbol
\begin{equation}
\label{ordre}
\alpha\leq\beta
\end{equation}
on $s$-uples of nonnegative integers will denote
the {\bf degree-antialphabetical \rm (total) \bf order}, defined by
$\alpha=(\alpha _1,\dots ,\alpha _s)<\beta=(\beta _1,\dots ,\beta _s)$
when
$\bigg( |\alpha |<|\beta |\bigg)$
or
$\bigg( |\alpha |=|\beta |$ and $\exists r\in\{ 1,\dots ,s\}$ such that
$\alpha _r<\beta _r$ and $\alpha _t=\beta _t$ for any $t>r\bigg)$.
For this order, $\delta _1<\delta _2<\dots <\delta _s
<2\delta _1<\delta _1+\delta _2\cdots$.

When $\alpha=(\alpha _1,\dots ,\alpha _s)$ is a $s$-uple of reals, 
the inequality
$$
\alpha \geq 0
$$
will mean that all the numbers $\alpha _k$ are $\geq 0$.

\masection{Main results}
\label{results-sec}

As it was briefly explained in Section~\ref{intro}, the method used to study
the asymptotics of a P\'olya process $(X_n)_n$ relies on estimates of its
moments in a Jordan basis, namely $E\uu ^\alpha (X_n)$,
$\alpha\in (\g Z_{\geq 0})^s$.
To this end, as it is developed in Subsection~\ref{phimoments-subsec},
it is natural to consider the transition operator $\Phi$ as it
was defined by Equation~(\ref{phi}).
Proposition~\ref{phiordre} is the first result on the action of $\Phi$
on polynomials.
One can find a proof of it in Subsection~\ref{phiOnPolynoms-subsec}.

\begin{Prop}
\label{phiordre}
For any choice of a Jordan basis $(u_k)_{1\leq k\leq s}$ of linear forms
of a P\'olya process and for every $\alpha \in (\g Z _{\geq 0})^s$,
$$
\Phi (\uu ^\alpha )-\langle\alpha ,\lambda\rangle \uu ^\alpha
\in \Vect\{ \uu ^\beta ,~\beta <\alpha\}.
$$
\end{Prop}

\noindent
The complex numbers $\langle\alpha ,\lambda\rangle$ were defined
in~(\ref{pdtscal}).
An immediate  consequence of this proposition is the $\Phi$-stability of the
finite-dimensional polynomial subspace
\begin{equation}
\label{defSalpha}
S_\alpha=\Vect\{\uu ^\beta ,~\beta\leq\alpha\}
\end{equation}
for any $\alpha \in (\g Z _{\geq 0})^s$.
These subspaces form an increasing sequence whose union
is the space $S(V)$ of all polynomial functions on $V$, so that
Proposition~\ref{phiordre} asserts that the eigenvalues of $\Phi$ on $S(V)$
are exactly all numbers
$\langle\alpha ,\lambda\rangle$, $\alpha \in (\g Z _{\geq 0})^s$
(in the (ordered) basis $(\uu _\beta )_{\beta\leq\alpha}$ of any $S_\alpha$,
the matrix of $\Phi$ is triangular).

\vskip 5pt
Notation: if $\Psi$ is an endomorphism of any vector space,
we will denote by $\ker\Psi ^\infty$ the characteristic space of $\Psi$
associated with zero, that is
\begin{equation}
\label{kerinfini}
\ker\Psi ^\infty =\bigcup _{p\geq 0}\ker\Psi ^p.
\end{equation}
We will use the notation $\Phi$ to refer to $\Phi$ itself as well as to
the endomorphism induced by $\Phi$ on $S(V)$ or on some stable subspace.
Decomposition of all $S_\alpha$ as direct sums of characteristic subspaces
of $\Phi$ leads to the splitting
$$
S(V)=\bigoplus _{z\in\g C}\ker (\Phi -z)^\infty.
$$
As it was announced in Section~\ref{intro}, we can now properly define the
reduced polynomials.

\begin{Def}
\label{defQalpha}
For any choice of a Jordan basis $(u_k)_{1\leq k\leq s}$ of linear forms of
a P\'olya process and
for any $\alpha \in (\g Z _{\geq 0})^s$,
the {\bf reduced polynomial of rank $\alpha$}
is the projection of $\uu ^\alpha$
on $\ker (\Phi -\langle\alpha ,\lambda\rangle )^\infty$
parallel to
$\bigoplus _{z\neq\langle\alpha ,\lambda\rangle}\ker (\Phi -z)^\infty$.
It will be denoted by $Q_\alpha$.
\end{Def}
Properties of reduced polynomial will be further developed in
Section~\ref{phi-sec}.
In particular, it will be explained how one can compute them inductively
(see~(\ref{inductionQalpha})).
They admit sometimes closed formulae
(see~(\ref{Qdelta1}), \cite{Barcelona} and~(\ref{triangle2})).
It follows from its definition that $Q_\alpha$ belongs to
$\ker (\Phi -\langle\alpha ,\lambda\rangle )^\infty$;
the number $\nu _\alpha$ defined just below is its index of nilpotence
in this characteristic space.
In particular, $\nu _\alpha =0$ if, and only if $Q_\alpha$ is eigenvector of
$\Phi$.
Proposition~\ref{nualpha} in Subection~\ref{nualpha-subsec} shows how one can
easily compute this number for any power of large projections.

\begin{Def}
\label{defnulpha-def}
For every $\alpha\in (\g Z _{\geq 0})^s$, the nonnegative integer $\nu _\alpha$
is defined by
\begin{equation}
\label{defnualpha}
\nu _\alpha=\max\{ p\geq 0,~(\Phi -\langle\alpha ,\lambda\rangle )^p(Q_\alpha )
\neq 0\} .
\end{equation}
\end{Def}

\vskip 5pt
These facts, definitions and notations being given, we claim the following
three main results of the article.

\begin{Th}{\bf (Joint moments of small or large projections)}
\label{jointmoments}
Let $(u_k)_{1\leq k\leq s}$ be a Jordan basis of linear forms of a P\'olya
process $(X_n)_n$.
Let $\alpha\in (\g Z_{\geq 0})^s$.

1- If $\alpha$ is a power of small projections, then there exists some
nonnegative integer $\nu$ such that
$$
E\uu ^\alpha (X_n)\in O\left( n^{|\alpha |/2}\log ^\nu n\right)
$$
as $n$ tends to infinity.

2- If $\alpha$ is a power of large projections, then there exists a complex
number $c$ such that
$$
E\uu ^\alpha \left( X_n\right)=
cn^{\langle \alpha ,\lambda \rangle}\log ^{\nu _\alpha}n
+o\left( n^{\Re\langle \alpha ,\lambda \rangle }\log ^{\nu _\alpha}n\right)
$$
as $n$ tends to infinity.

3- If $\alpha$ is a semisimple power of large projections, then
$$
E\uu ^\alpha \left( X_n\right) =
n^{\langle \alpha ,\lambda \rangle}
\frac{\Gamma (\tau _1)}{\Gamma (\tau _1+\langle \alpha ,\lambda \rangle )}
Q_\alpha (X_1)
+o\left( n^{\Re\langle \alpha ,\lambda \rangle}\right)
$$
as $n$ tends to infinity, where $Q_\alpha$ is the reduced polynomial of rank
$\alpha$ relative to the Jordan basis $(u_k)_{1\leq k\leq s}$.
\end{Th}
Constant $c$ in Assertion~{\it 2-} has an explicit form given in
Remark~\ref{constantc}.
The proof of Theorem~\ref{jointmoments} can be found in
Section~\ref{jointmoments-sec}.
It is based on a careful study of coordinates of the $\uu$-monomials in the
basis of reduced polynomials, which is developed in
Subsections~\ref{cones-subsec} and~\ref{continued-subsec}.

Although it is not formally necessary, we give two different statements on
the asymptotics of large P\'olya processes, respectively when the process
is principally semisimple or not.
Their proofs can be found in Section~\ref{main-sec}.
They are based on Theorem~\ref{jointmoments} and use martingale techniques
(quadratic variation, Burkholder Inequality).

\begin{Th}
{\bf (Asymptotics of large and principally semisimple P\'olya processes)}
\label{lpss}
Suppose that a P\'olya process $(X_n)_n$ is large and principally semisimple.
Fix a Jordan basis $(u_k)_{1\leq k\leq s}$ of linear forms such that
$u_1,\dots ,u_r$ ($2\leq r\leq s$) are all the eigenforms of the basis
that are associated with roots\footnote{
In short, if $1$ is multiple root, $\lambda _1=\dots =\lambda _r=1$;
otherwise, $\frac 12<\Re\lambda _2=\dots =\Re\lambda _r=\sigma _2<1$.
See~(\ref{defsigma2}), definition of $\sigma _2$.
}
$\lambda _1=1,\lambda _2, \dots ,\lambda _r$ having
a real part $\geq\sigma _2$.

Then, with notations~(\ref{vk}) and~(\ref{lambda}) of
Section~\ref{Preliminaries},
there exist unique (complex-valued) random variables $W_2,\dots ,W_r$ such that
\begin{equation}
\label{asymptlpss}
\displaystyle
X_n=nv_1+\sum _{2\leq k\leq r}n^{\lambda _k}W_kv_k
+o\left( n^{\sigma _2}\right),
\end{equation}
the small $o$ being almost sure and in ${\rm L}^p$ for every $p\geq 1$.
Furthermore, if one denotes by $(Q_\alpha )_{\alpha\in (\g Z_{\geq 0})^s}$
the reduced polynomials relative to the Jordan basis $(u_k)_{k}$,
all joint moments of the random variables $W_2,\dots ,W_r$
exist and are given by: for all
$\alpha _2,\dots ,\alpha _r\in \g Z_{\geq 0}$,
$$
E\left( \prod _{2\leq k\leq r}W_k^{\alpha _k}\right)
=\frac{\Gamma (\tau _1)}{\Gamma (\tau _1+\langle \alpha ,\lambda \rangle )}
{Q_\alpha (X_1)}
$$
where $\alpha =\sum _{2\leq k\leq r}\alpha _k\delta _k
=(0,\alpha _2,\dots ,\alpha _r,0,\dots )$.
\end{Th}

\begin{Th}
{\bf (Asymptotics of large and principally nonsemisimple P\'olya processes)}
\label{lpnss}
Suppose that the P\'olya process $(X_n)_n$ is large and principally
nonsemisimple.
Fix a Jordan basis $(u_k)_{1\leq k\leq s}$ of linear forms;
let $J_2,\dots ,J_r$ be the principal blocks of
indices\footnote{
In other words, if $J$ is any Jordan block of $A$ in the $u_k$'s basis,
$J$ is $1$ or one of the $J_k$'s, or the size of $J$ is $\leq\nu$, or the
root of $J$ has a real part $<\sigma _2$.
See Definition \ref{defPrincipalBlock} (principal blocks).
}
and $\nu +1$ the common size of the $J_k$'s ($\nu\geq 1$).

Then, with notations (\ref{vk}) and (\ref{lambda}) of
Section~\ref{Preliminaries},
there exist unique (complex-valued) random variables $W_2,\dots ,W_r$ such that
\begin{equation}
\label{asymptlpnss}
\displaystyle
X_n=nv_1+\frac 1{\nu !}
\log ^\nu n\sum _{2\leq k\leq r}n^{\lambda (J_k)}W_kv_{\max J_k}
+o\left( n^{\sigma _2}\log ^\nu n\right),
\end{equation}
the small $o$ being almost sure and in ${\rm L}^p$ for every $p\geq 1$.
Furthermore, if one denotes by $(Q_\alpha )_{\alpha\in (\g Z_{\geq 0})^s}$
the reduced polynomials relative to the Jordan basis $(u_k)_k$,
all joint moments of the random variables $W_2,\dots ,W_r$
exist and are given by: for all
$\alpha _2,\dots ,\alpha _r\in \g Z_{\geq 0}$,
$$
E\left( \prod _{2\leq k\leq r}W_k^{\alpha _k}\right)
=\frac{\Gamma (\tau _1)}{\Gamma (\tau _1+\langle\alpha ,\lambda\rangle )}
Q_\alpha (X_1)
$$
where
$\alpha =\sum _{2\leq k\leq r}\alpha _k\delta _{\min J_k}$.
\end{Th}

\masection{Transition operator}
\label{phi-sec}

Let $(X_n)_n$ be a P\'olya process given by its increment vectors
$(w_k)_{1\leq k\leq s}$ and its basis of linear forms $(l_k)_{1\leq k\leq s}$
submitted to hypotheses of Definition~\ref{defiPolya}.
We recall here the definition of its associated transition operator $\Phi$
as it was given in Section~\ref{intro}:
if $f:V\to W$ is any $W$-valued function where $W$ is any real vector space,
$\forall v\in V$,
\begin{equation*}
\Phi (f)(v)=
\sum _{1\leq k\leq s}l_k(v)\bigg[ f(v+w_k)-f(v)\bigg].
\end{equation*}

\subsection{Transition operator $\Phi$ and computation of moments}
\label{phimoments-subsec}

Proposition \ref{esper} expresses the expectation of any $f(X_n)$ in
terms of $f$, of iterations of the transition operator $\Phi$
and of $X_1$, initial value of the process.
Polynomials $\gamma _{\tau _1,n}$ with rational coefficients and one
variable were defined by Equation~(\ref{gamma}).

\begin{Prop}
\label{esper}
If $f:V\to W$ is any measurable function
taking values in some real (or complex) vector space $W$,
then for all $n\geq 1$,
\begin{equation}
\label{esperfXn}
Ef(X_n)=\gamma _{\tau _1,n}(\Phi)(f)(X_1).
\end{equation}
\end{Prop}

\begin{pff}
It follows immediately from~(\ref{probaTransition}) that the expectation of
$f(X_{n+1})$ conditionally to the state at time $n$ is
\begin{eqnarray*}
E^{{\cal F}_n}f(X_{n+1})&=
&\sum _{1\leq k\leq s}\frac1{n+\tau _1-1}l_k(X_n)f(X_n+w_k)\\
&=&f(X_n)+
\frac1{n+\tau _1-1}
\sum _{1\leq k\leq s}l_k(X_n)\bigg(f(X_n+w_k)-f(X_n)\bigg).
\end{eqnarray*}
By definition of the transition operator $\Phi$, this formula can be written as
\begin{equation}
\label{espeCondi}
E^{{\cal F}_n}f(X_{n+1})=\left( \Id +\frac1{n+\tau _1-1}\Phi\right)(f)(X_n);
\end{equation}
taking the expectation leads to the result after a straightforward induction.
\end{pff}\QED

\vskip 5pt
It follows from Proposition \ref{esper} that the asymptotic weak behaviour
of the process, or at least the asymptotic behaviour of its
moments is reachable by decompositions
of the operator $\Phi$ on suitable function spaces.
Corollary \ref{fjordan} is the first step in this direction, stating the
result for functions that belong to finite dimensional stable subpaces.

\begin{Cor}
\label{fjordan}
Let $f:V\to W$ be a measurable function taking values in some real (or complex)
vector space $W$.

1- If $f$ is an eigenfunction of $\Phi$ associated with the
(real or complex) eigenvalue~$z$, that is if
$\Phi (f)=z f$, then
$$
Ef(X_n)
=n^z\frac{\Gamma (\tau _1)}{\Gamma (\tau _1+z)}f(X_1)+O\left(n^{z -1}\right)
$$
as $n$ tends to infinity ($\Gamma$ is Euler's function).

2- Assume that $f$ is nonzero and belongs to some $\Phi$-stable subspace
$\cal S$ of measurable functions $V\to W$ and that the operator induced by
$\Phi$ on $\cal S$ is a sum $z \Id _{\cal S} +\Phi _N$, where $\Phi _N$ is a
nonzero nilpotent operator on $\cal S$ and $z$ a complex number.
Let $\nu$ be the positive integer such that $\Phi _N^\nu (f)\neq 0$ and
$\Phi _N^{\nu +1}(f)=0$.
Then,
$$
Ef(X_n)=
\frac{n^z\log ^\nu n}{\nu !}\frac{\Gamma (\tau _1)}{\Gamma (\tau _1+z)}
\Phi _N^\nu (f)(X_1)
+O\left( n^z\log ^{\nu -1}n\right)
$$
as $n$ tends to infinity.
\end{Cor}

\begin{pff}
{\it 1-} It follows from Proposition \ref{esper} that
$Ef(X_n)=\gamma _{\tau _1,n}(z )\times f(X_1)$.
Note that, as soon as the terms are defined,
\begin{equation}
\label{gammaStirling}
\gamma _{\tau _1,n}(t )=
\frac {\Gamma (\tau _1)}{\Gamma (\tau _1+t)}
\frac{\Gamma (n+\tau _1-1+t)}{\Gamma (n+\tau _1-1)},
\end{equation}
so that the result is a consequence of Stirling Formula.

\noindent
{\it 2-}
Taylor expansion of $\gamma _{\tau _1,n}(z\Id +\Phi _N)$ leads to
$$
Ef(X_n)=\sum _{p\geq 0}\frac 1{p!}\gamma _{\tau _1,n}^{(p)}(z )
\Phi _N^p(f)(X_1)
$$
(finite sum), where $\gamma _{\tau _1,n}^{(p)}$ denotes the $p$-th
derivative of $\gamma _{\tau _1,n}$.
Besides, if $p$ is any positive integer,
\begin{equation}
\label{derivGamma}
\gamma _{\tau _1,n}^{(p)}(z)=
n^z\log ^pn\frac{\Gamma (\tau _1)}{\Gamma (\tau _1+z)}
+O\left( n^z\log ^{p-1}n\right)
\end{equation}
when $n$ tends to infinity,
as can be shown by Stirling formula (see (\ref{gammaStirling}))
and an elementary induction starting from
the computation of $\gamma _{\tau _1,n}$'s logarithmic derivative.
These two facts imply the result.
\end{pff}\QED

\begin{rem}
As it is written, Corollary \ref{fjordan} is valid only if the complex number
$\tau _1+z$ is not a nonpositive integer.
We adopt the convention $1/\Gamma (w)=0$
when $w\in\g Z_{\leq 0}$, so that this corollary is valid in all cases.
\end{rem}

\begin{rem}
\label{flinear}
If $f:V\to W$ is linear, formula (\ref{phi}) implies that
$\Phi (f)=f\circ A$.
In that particular case, formula (\ref{esperfXn}) gives
$Ef(X_n)=f\circ\gamma _{\tau _1,n}(A)(X_1)$.
This fact will be used in the proofs of Theorems~\ref{lpss} an~\ref{lpnss}
when $f$ is a linear combination of projections $\pi _k$
(see Section~\ref{main-sec}).
\end{rem}

\subsection{Action of $\Phi$ on polynomials}
\label{phiOnPolynoms-subsec}

Because of Condition~(\ref{hypobalance}) in the definition of a P\'olya
process, none of the vectors $w_k$ is zero.
For any $k$, if $f$ is a function defined on $V$, we denote by
$\partial f/\partial w_k$,
when it exists, the derivative of $f$ along the direction
carried by the vector $w_k$.
With this notation, we associate with the finite difference operator $\Phi$
the {\bf differential operator} $\Phi _\partial$ defined by
\begin{equation}
\label{defPhipartial}
\Phi _\partial (f)(v)=\sum _{1\leq k\leq s}l_k(v)
\frac {\partial f}{\partial w_k}(v)
\end{equation}
for every function $f$ defined on $V$ and derivable at each point
along the directions carried by the vectors $w_k$'s.
When $f$ is differentiable, $\Phi _\partial (f)$ can be viewed as a
``first approximation'' of $\Phi (f)$.
As derivation behaves good with respect to product of functions when finite
differentiation does not, $\Phi _\partial (f)$ is helpful for the understanding
of $\Phi$'s  action on polynomials.

\begin{rem}
\label{differential}
The differential operator can be written as $\Phi _\partial (f)(v)=\D f_v.Av$
for any differentiable function $f$, where $\D f_v$ denotes the differential
of $f$ at point $v$.
This can be readily seen from the formula
$\D f_v.w_k=\frac{\partial f}{\partial w_k}(v)$.
\end{rem}

\begin{Prop}
\label{phidiffpropre}
{\bf (Action of $\Phi _\partial$ on the $\uu$-monomials)}
For any choice of a Jordan basis $(u_k)_{1\leq k\leq s}$ of linear forms
of a P\'olya process,

1- for every $\alpha \in (\g Z _{\geq 0})^s$,
$$
\Phi _\partial (\uu ^\alpha )-\langle\alpha ,\lambda\rangle \uu ^\alpha
\in \Vect\{ \uu ^\beta ,\beta <\alpha\} ;
$$

2- if $\alpha \in (\g Z _{\geq 0})^s$ is a semisimple power,
then
$\Phi _\partial (\uu ^\alpha )=\langle\alpha ,\lambda\rangle\uu ^\alpha$.
\end{Prop}

\begin{pff}
$\Phi _\partial$ is a derivation, as can be seen directly or from
Remark~\ref{differential}.
In particular, for any $\alpha\in (\g Z _{\geq 0})^s$,
\begin{equation}
\label{diffualpha}
\Phi _\partial (\uu ^\alpha )
=\sum _{k=1}^s\alpha _k\uu ^{\alpha -\delta _k}\Phi _\partial (u_k).
\end{equation}
Besides, as any $u_k$ is linear, $\Phi _\partial(u_k)=u_k\circ A$.
The conclusion follows from Jordan basis' Definition~\ref{jordanBasis-def}
(the degree-antialphabetical order on $s$-uples is defined in (\ref{ordre})
at the end of Subsection \ref{powers}).
\end{pff}\QED

\begin{rem}
One can formally extend the result of {\it 2-} in
Proposition~\ref{phidiffpropre} to any family of {\it complex} numbers
$\alpha _1, \dots, \alpha _s$ when $\forall k$,
$k\neq 0\Longrightarrow u_k$ \it is eigenform of $A$.
\rm This gives other eigenfunctions of $\Phi _\partial$, defined on suitable
open subsets of $V$ or $V_{\g C}$ (usual topology).
\end{rem}

We can now prove Proposition~\ref{phiordre}, as it was announced in
Section~\ref{results-sec}.
It appears as a direct consequence of Proposition~\ref{phidiffpropre}.

\vskip 5pt\noindent
{\sc Proof of Proposition~\ref{phiordre}.}\ 
\label{preuvephiordre}
The family $(\uu ^\beta ) _{|\beta |\leq |\alpha |-1}$ constitutes a
basis of polynomials of degree $\leq |\alpha |-1$.
Hence, if $F$ denotes the subspace $F=\Vect\{\uu ^\beta ,~\beta <\alpha\}$,
Taylor formula implies that
$(\Phi -\Phi _\partial )(\uu ^\alpha )\in F$.
Moreover,
$\Phi _\partial (\uu ^\alpha )-\langle\alpha ,\lambda\rangle \uu ^\alpha\in F$
because of Proposition \ref{phidiffpropre}.
This completes the proof.
\QED

\subsection{Reduced polynomials}
\label{Qalpha-subsec}

Choose a Jordan basis of linear forms $(u_k)_{1\leq k\leq s}$ of a P\'olya
process.
For any $\alpha \in (\g Z _{\geq 0})^s$, the reduced polynomial of rank
$\alpha$, denoted by $Q_\alpha$, was
defined in Definition~\ref{defQalpha} as the projection of $\uu ^\alpha$
on $\ker (\Phi -\langle\alpha ,\lambda\rangle )^\infty$
parallel to
$\bigoplus _{z\neq\langle\alpha ,\lambda\rangle}\ker (\Phi -z)^\infty$
(see~(\ref{kerinfini}) for the meaning of notation $\ker\psi ^\infty$).
Properties of these polynomials that are listed in
Proposition~\ref{proprietesQalpha} will be used in the sequel.
Subspaces $S_\alpha$ were defined in~(\ref{defSalpha}).

\begin{Prop}
\label{proprietesQalpha}
Let $\alpha \in (\g Z _{\geq 0})^s$.

\item (1)
$Q_0=1$ and $Q_\alpha=\uu ^\alpha$ if $|\alpha |=1$;

\item (2)
$\{Q_\beta,~\beta\leq\alpha\}$ is a basis of $S_\alpha$;

\item (3)
for every $z\in\g C$, $\{Q_\alpha ,~\langle\alpha ,\lambda\rangle =z\}$
is a basis of $\ker (\Phi -z)^\infty$;

\item (4)
$Q_\alpha-\uu ^\alpha\in\Vect\{ Q_\beta,~\beta <\alpha ,
~\langle\beta ,\lambda\rangle\neq\langle\alpha ,\lambda\rangle\}$;

\item (5)
$\Phi (Q_\alpha )-\langle\alpha ,\lambda \rangle Q_\alpha\in
\Vect\{Q_\beta ,~\beta <\alpha ,
~\langle\beta ,\lambda \rangle =\langle\alpha ,\lambda \rangle\}$.
\end{Prop}

\begin{pff}
{\it (1)} comes directly from the definition of a Jordan basis and
{\it (2)} from the $\Phi$-stability of subspaces $S_\alpha$
(see~(\ref{defSalpha})).
Any $Q_\alpha$ belongs to the characteristic space
$\ker (\Phi -\langle\alpha ,\lambda\rangle )^\infty$ and the eigenvalues
of the restriction of $\Phi$ on polynomials are exactly the
$\langle\alpha ,\lambda\rangle$ (see Section~\ref{results-sec}, consequences
of Proposition~\ref{phiordre}).
These facts imply {\it (3)}.
Property {\it (4)} is obvious from $Q_\alpha$'s definition, when {\it (5)}
follows from {\it (3)} and Proposition~\ref{phiordre}.
\end{pff}
\QED

Assertions {\it (4)} and {\it (5)} in Proposition~\ref{proprietesQalpha}
can be used to compute the reduced
polynomials inductively (see Remark~\ref{Qalphainduction} below).
Let's define, as it was announced in Section~\ref{intro},
the complex numbers $q_{\alpha ,\beta}$ by the relations
$\uu ^\alpha =Q_\alpha
+\sum _{\beta <\alpha }q_{\alpha ,\beta}Q_\beta$,
their existence and unicity being guaranteed by Assertions~{\it (2)}
and~{\it (4)} in Proposition~\ref{proprietesQalpha}.
Moreover, because of {\it (4)}, $q_{\alpha ,\beta}=0$ as soon as
$\langle\beta ,\lambda\rangle =\langle\alpha ,\lambda\rangle$, so that
\begin{equation}
\label{qalphabeta}
\uu ^\alpha =Q_\alpha
+\sum _
{\beta <\alpha ,~\langle\beta ,\lambda\rangle\neq\langle\alpha ,\lambda\rangle}
q_{\alpha ,\beta}Q_\beta .
\end{equation}
This relation, still too rough to lead to the main results on asymptotics of
large P\'olya processes, will be refined in Subsection~\ref{continued-subsec}.

We end the present subsection by giving two remarks concerning the inductive
computation of all reduced polynomials and a closed form for projections of the
powers of $u_1$.

\begin{rem}
\label{Qalphainduction}
{\bf Inductive computation of $Q_{\alpha}$'s.}

\noindent
In the general case, the numbers $q_{\alpha ,\beta}$
and the numbers $p_{\alpha ,\beta}$ defined by
$$
(\Phi -\langle\alpha ,\lambda\rangle )(Q_\alpha )
=\sum _{\beta <\alpha}p_{\alpha ,\beta}Q_\beta
=\sum _{\beta <\alpha ,
~\langle\beta ,\lambda\rangle =\langle\alpha ,\lambda\rangle}
p_{\alpha ,\beta}Q_\beta
$$
(see {\it (5)} in Proposition~\ref{proprietesQalpha})
can be inductively computed (and implemented) the
following way.
We denote by $r_{\alpha ,\beta}$ the complex numbers defined by
\begin{equation}
\label{ralphabeta}
(\Phi -\langle\alpha ,\lambda\rangle )(\uu ^\alpha )
=\sum _{\beta <\alpha}r_{\alpha ,\beta}Q_\beta,
\end{equation}
that can be deduced by plain computation of
$(\Phi -\langle\alpha ,\lambda\rangle )(\uu ^\alpha )$ and its expansion
in the $(Q_\beta )_{\beta <\alpha}$ basis with the help of formula
(\ref{qalphabeta}), the corresponding numbers $q_{\beta ,\gamma}$
being known by induction.
Write two expressions of $(\Phi -\langle\alpha ,\lambda\rangle )(\uu ^\alpha )$
with formulae (\ref{qalphabeta}) and (\ref{ralphabeta}) and identify
the coordinates in the $(Q_\beta )$ basis.
This provides the following equations with $p_{\alpha ,\beta}$ and
$q_{\alpha ,\beta}$ as unknowns:
\begin{equation}
\label{inductionQalpha}
\left\{
\begin{array}{rcl}
\langle\beta ,\lambda\rangle =\langle\alpha ,\lambda\rangle &
\Longrightarrow&
r_{\alpha ,\beta}=p_{\alpha ,\beta}\\
\langle\beta ,\lambda\rangle\neq\langle\alpha ,\lambda\rangle &
\Longrightarrow&
r_{\alpha ,\beta}=\left(
\langle\beta ,\lambda\rangle -\langle\alpha ,\lambda\rangle\right)
q_{\alpha ,\beta}\displaystyle
+\sum _{\beta <\gamma <\alpha}q_{\alpha ,\gamma}p_{\gamma ,\beta}.\\
\end{array}
\right.
\end{equation}
The expansion of $Q_\alpha$ in the $(\uu ^\beta)$ basis can
be obtained by reversing the triangular system written in (\ref{qalphabeta}).
All these computations can be handled by means of symbolic computation.
\end{rem}

\begin{rem}
\label{Qdelta1-rem}
{\bf Closed formula for $Q_{p\delta _1}$'s.}

\noindent
An immediate computation shows that the reduced polynomials corresponding
to powers of $u_1$ are the same ones for all P\'olya processes:
for any integer $p\geq 0$,
$\Phi \left[ u_1(u_1+1)\dots (u_1+p-1)\right] =pu_1(u_1+1)\dots (u_1+p-1)$,
so that if follows from Proposition~\ref{proprietesQalpha} that
$Q_{p\delta _1}=u_1(u_1+1)\dots (u_1+p-1)$.
The powers of $u_1$ are thus always expressed in terms of reduced polynomials
$Q_{p\delta _1}$'s by means of Stirling numbers of the second kind
(for this inversion formula, see {\it e.g.} \cite{ConcreteMath}):
\begin{equation}
\label{Qdelta1}
u_1^p=\sum _{k=1}^p(-1)^{p-k}\stirling2 pkQ_{k\delta _1}.
\end{equation}
This common formula has to be related to the non random drift, consequence
of~(\ref{balance}): $\forall n,~u_1(X_n)=n+\tau _1-1$.

\end{rem}

\subsection{Cones in the space of powers}
\label{cones-subsec}

As it was explained in Section~\ref{intro}, it follows from reduced
polynomial's definition that the behaviour of the $EQ_\beta (X_n)$ are ruled
by Corollary~\ref{fjordan}.
Thus, the asymptotics of the $\uu$-moment
\begin{equation}
\label{espUespQ}
E\uu ^\alpha (X_n)=EQ_\alpha (X_n)+
\sum _
{\beta <\alpha ,~\langle\beta ,\lambda\rangle\neq\langle\alpha ,\lambda\rangle}
q_{\alpha ,\beta}EQ_\beta (X_n),
\end{equation}
when $n$ goes off to infinity,
depends on the answer to the following two questions:
\begin{center}
\begin{minipage}{11,5 truecm}
\item (1)
which $q_{\alpha ,\beta}$ are zero in Relation (\ref{qalphabeta})?

\item (2)
For a given $\alpha$, which $\Re\langle\beta ,\lambda\rangle$ is maximal
among indices $\beta <\alpha$ such that $q_{\alpha ,\beta}\neq 0$?
\end{minipage}
\end{center}
Optimal answer for the most general P\'olya process is expressed in terms of
a rational cone $\Sigma$ and a rational polyhedron $A_\alpha$ in the ``space
of powers'' $\g R^s=\g Z^s\otimes\g R$.
The two following paragraphs are devoted to these subsets;
at the end of each of them, we give properties of the number
$\langle\beta ,\lambda\rangle$ when $\beta$ belongs respectively to
$\Sigma$ or some $A_\alpha$
(Propositions~\ref{majorationsRe} and~\ref{propAalphalambda}).

Note, as suggested by the anonymous referee, that the argument given in
Remark~\ref{GammaReferee} enables one to by-pass Sections~\ref{Sigma-subsubsec}
and~\ref{continued-subsec} giving the construction of the cone $\Sigma$ and
the study of its properties.

\subsubsection{Cone $\Sigma$}
\label{Sigma-subsubsec}

{\it Notations:}
if $I\subseteq\{1,\dots ,s\}$ and
$(i,j)\in\{1,\dots ,s\}^2$, we adopt the notations
\begin{eqnarray*}
\delta _I=\sum_{1\leq i\leq s}\delta _i+\sum _{i\in I}\delta _i\in\g R^s
&\rm and&
\delta _I^*=\sum_{1\leq i\leq s}dx_i+\sum _{i\in I}dx_i\in{\g R^s}^*,\\
\\
\delta _{(i,j)}=2\delta _i-\delta _j\in\g R^s&\rm and&
\delta _{(i,j)}^*=2dx_i-dx_j\in{\g R^s}^*
\end{eqnarray*}
where $dx_k$ denotes the $k$-th coordinate form $(x_1,\dots ,x_s)\mapsto x_k$
in the dual space $\g R^{s*}$ and where $\delta _k$ is the $k$-th vector of
the canonical basis of $\g R^s$, already defined in~(\ref{deltak}).

\begin{Def}
\label{Sigma}
We denote by $\Sigma$
the polyhedral cone of $\g R^s$ spanned by the
$s(s-1)$
vectors $\delta _{(i,j)}$
for all ordered pairs $(i,j)$ of distinct elements, {\it i.e.}
\begin{equation}\label{defSigma}
\Sigma =\sum_ {(i,j)\in\{1,\dots s\},~i\neq j} \g R_{\geq 0}\delta _{(i,j)}.
\end{equation}
\end{Def}
This cone is convex, and the half-lines spanned by vectors $\delta _{(i,j)}$
are extremal (edges).
As usual, we define the dual cone $\check\Sigma$ of $\Sigma$ as
$$
\check\Sigma =\{x\in \g R^s, \forall y\in\Sigma , \langle x,y\rangle\geq 0\},
$$
identified to the cone of all linear forms on $\g R^s$ that are nonnegative
on $\Sigma$, via the bijective linear application
$x\in\g R^s\mapsto \langle x,.\rangle\in{\g R^s}^*$
(the symbol $\langle x,y\rangle$ denotes the standard scalar product of $x$ and
$y$ in $\g R^s$).
Lemma \ref{dual} describes the dual cone $\check\Sigma$ as a minimal
intersection of hyperplanes (faces) and gives a system of minimal generators
(edges).
Corollary \ref{facesSigma} just transcribes Lemma \ref{dual} in the
$\Sigma$-side and gives the equations of the faces of $\Sigma$.
We give a complete geometrical description of $\Sigma$;
it presents some ``universal'' character, as shown in Remark~\ref{equivSigma}.

\begin{Lem}
\label{dual}
{\bf (faces and edges of $\check\Sigma$)}
$$
\check\Sigma =\displaystyle\bigcap _{(i,j)\in\{1,\dots s\},~i\neq j}
\{x\in\g R^s,~\delta _{(i,j)}^*(x)\geq 0\}
=
\summ {I\subseteq \{1,\dots ,s\}}{1\leq\#I\leq s-1}
\g R_{\geq 0}\delta _I.
$$
\end{Lem}

\begin{pff}
The first equality that describes the faces of $\check\Sigma$ comes directly
from (\ref{defSigma}).
For every permutation $w\in{\mathfrak S}_s$, let $\tau _w$ be the simplicial
cone defined by
$$
\tau _w=
\{x\in\g R^s,~x_{w(s)}\leq x_{w(s-1)}\leq\dots\leq x_{w(1)}\leq 2x_{w(s)}\}.
$$
The cones $\tau _w$ provide a subdivision of $\check\Sigma$ in $s!$ simplicial
cones -this subdivision is the intersection of $\check\Sigma$ with the
barycentric subdivision of the first quadrant of $\g R^s$.
Each $\tau _w$ is the image of $\tau _1=\tau _{\Id}$ by the permutation of
coordinates induced by $w$
(and $\tau _1$ is a fundamental domain for the group action of
${\mathfrak S}_s$ on $\check\Sigma$ by permutations of coordinates).
Because of the elementary computation
\begin{eqnarray*}
(x_1,x_2,x_3,x_4)=&(2x_4-x_1)(1,1,1,1)\\
&+(x_1-x_2)(2,1,1,1)\\
&+(x_2-x_3)(2,2,1,1)\\
&+(x_3-x_4)(2,2,2,1)
\end{eqnarray*}
that can be straightforwardly generalized in all dimensions, one sees that the
edges of $\tau _1$ are spanned by $(1,\dots ,1)=\delta _{\emptyset}$
and $\delta _{\{1\}}$, $\delta _{\{1,2\}}$,\dots ,
$\delta _{\{1,\dots, s-1\}}$.
The images of these last $s-1$ vectors under permutations of coordinates are
exactly the $\delta _I$, where $I\neq\emptyset$ and $I\neq\{1,\dots ,s\}$.
This completes the proof.
\end{pff}\QED

\begin{rem}
\label{deltavide}
Vectors $\delta _\emptyset =\sum _{1\leq k\leq s}\delta _k$ and
$\delta _{\{ 1,\dots ,s\}}=2\delta _\emptyset$ belong to $\check\Sigma$,
and this fact will be used in the sequel.
They do not appear in the second sum of Lemma~\ref{dual} because they do not
span an edge of $\check\Sigma$.
On the contrary, the vector $\delta _I$ spans an edge of $\check\Sigma$ when
$I$ is neither empty nor the whole $\{ 1,\dots ,s\}$.
\end{rem}

\begin{Cor}
\label{facesSigma}
{\bf (faces of $\Sigma$)}
The cone $\Sigma$ has $2^s-2$ faces of dimension $s-1$, described as
\begin{equation}
\label{eqFacesSigma}
\Sigma =\bigcapp {I\subseteq \{1,\dots ,s\}}{1\leq\#I\leq s-1}
\{x\in\g R^s,~\delta _I^*(x)\geq 0\}.
\end{equation}
\end{Cor}

\begin{figure}[htbp]
\psfrag{s}{$\Sigma$}\psfrag{a}{$\alpha$}
\psfrag{a-s}{$\alpha -\Sigma$}\psfrag{ss}{$\check\Sigma$}
\psfrag{e}{$\{\alpha -\eta,~\eta\geq 0,~|\eta |\geq 2\}$}
\psfrag{ee}{$\{\alpha -\eta +\delta _k,~\eta\geq 0,
~|\eta |\geq 2,~k\in\{1,2\}\}$}
\centerline{\includegraphics[width=20mm,height=20mm]{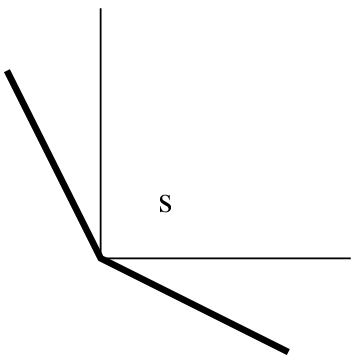}
\hskip 100pt
\includegraphics[width=15mm,height=15mm]{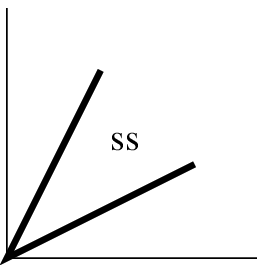}}
\vskip 15pt
\centerline{\hskip -30pt
\includegraphics[width=35mm,height=25mm]{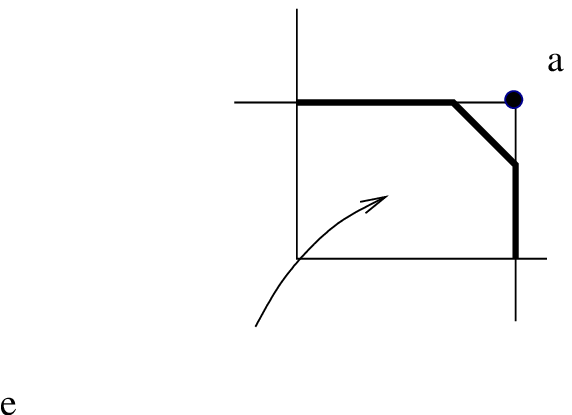}
\hskip 30pt
\includegraphics[width=40mm,height=30mm]{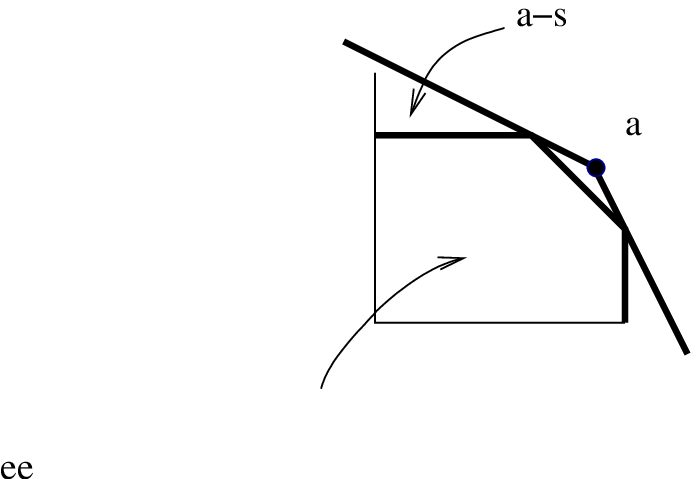}}
\caption{cone $\Sigma$ and related sets in dimension $2$}\label{Fi:sigma2}
\end{figure}

\begin{figure}[htbp]
\psfrag{1}{$\delta _1$}\psfrag{2}{$\delta _2$}\psfrag{3}{$\delta _3$}
\psfrag{4}{$(0,2,-1)$}\psfrag{5}{$(2,0,-1)$}\psfrag{6}{\hskip 5pt$(2,-1,0)$}
\psfrag{7}{\hskip -20pt $(0,-1,2)$}\psfrag{8}{$(-1,0,2)$}\psfrag{9}{$(-1,2,0)$}
\centerline{\includegraphics[width=50mm,height=40mm]{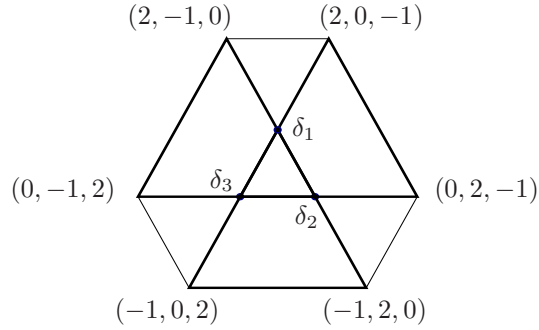}}
\caption{trace of $\Sigma$ (convex hull) and of
$\{\eta -\delta _k,~\eta\geq 0,~|\eta |\geq 2,~k\in\{1,2,3\}\}$
(union of three triangles)
on the hyperplane
$\{x_1+x_2+x_3=1\}$ of $\g R^3$}
\label{Fi:sigma3}
\end{figure}

In dimension two, $\Sigma$ is spanned by
$(2,-1)$, and $(-1,2)$
and $\check\Sigma$ by the forms
$2dx_1+dx_2$ and $dx_1+2dx_2$.
In dimension three, $\Sigma$ is spanned by
$(2,-1,0)$, $(-1,2,0)$, $(2,0,-1)$, $(-1,0,2)$, $(0,2,-1)$ and $(0,-1,2)$ and
the coordinates of the spanning forms of $\check\Sigma$ are
$(2,1,1)$, $(1,2,1)$, $(1,1,2)$, $(2,2,1)$, $(2,1,2)$ and $(1,2,2)$
in the canonical basis $(dx_1,dx_2,dx_3)$.
The numbers of edges of $\Sigma$ and $\check\Sigma$
coincide only in dimensions $2$ and $3$.
Figures \ref{Fi:sigma2} and \ref{Fi:sigma3} give pictures of $\Sigma$ in
dimensions $2$ and $3$;
in these figures, the comments that contain occurences of the greek letter
$\eta$ refer to further developments
(see Remark~\ref{remarkSigma} in Subsection~\ref{continued-subsec}).

\vskip 5pt\noindent
{\it Notation}: for any $B\subseteq\g R^s$ and $\alpha\in\g R^s$, we denote
\begin{equation}
\label{B-Sigma}
B-\Sigma =\{ B-\sigma,~\sigma\in\Sigma\}
{\rm ~~and~~}
\alpha -\Sigma =\{\alpha\}-\Sigma.
\end{equation}

\begin{Prop}
\label{majorationsRe}
Take a P\'olya process, choose any Jordan basis and denote
$\lambda =(\lambda _1,\dots ,\lambda _s)$ its root's $s$-uple
(see~(\ref{lambda})).
Let $\alpha\in (\g Z_{\geq 0})^s$ and
$\beta\in (\g Z_{\geq 0})^s\cap (\alpha -\Sigma )$.

1- If $\alpha$ is a power of large projections, then
$\alpha =\beta$ or
$\Re\langle\beta ,\lambda\rangle <\Re\langle\alpha ,\lambda\rangle$.

2- If $\alpha$ is a power of small projections, then
$\Re\langle\beta ,\lambda\rangle\leq \frac 12|\alpha |$.
\end{Prop}

\begin{pff}
See Definition~\ref{largesmall} in Subsection~\ref{powers} for definitions
of powers of large or small projections.
We denote $\sigma =\alpha -\beta\in\Sigma$ and split $\{ 1,\dots ,s\}$
into the three disjoint subsets
$I=\{k,~\Re (\lambda _k)\geq \frac 12,~\sigma _k\leq 0\}$,
$J=\{k,~\Re (\lambda _k)\geq \frac 12,~\sigma _k>0\}$
and
$K=\{k,~\Re (\lambda _k)<\frac 12\}$.

\noindent
{\it 1-}
If $\alpha$ is a power of large projections, then any $\alpha _k$ $(k\in K)$
vanishes.
Since $\alpha -\sigma\geq 0$, this implies that $\sigma _k\leq 0$ for all
$k\in K$.
Thus,
\begin{equation}
\label{inegGdeProj}
\Re\langle\sigma ,\lambda\rangle\geq
\sum _{k\in I}\sigma _k
+\frac 12\sum _{k\in J}\sigma _k
+\frac 12\sum _{k\in K}\sigma _k
=\frac 12\delta _I^*(\sigma ).
\end{equation}
Since $\sigma$ lies in $\Sigma$, the number $\delta _I^*(\sigma )$
is nonnegative (see~(\ref{eqFacesSigma})).
Besides, $\sigma _k\leq 0$ for any $k\in K$ so that
there exists some $k\in I\cup J$ such that $\sigma _k>0$
because the only point of $\Sigma$ with only nonpositive coordinates is $0$,
as can be seen on $\Sigma$'s equations (Corollary~\ref{facesSigma}).
Thus $J\neq\emptyset$ and the inequality of (\ref{inegGdeProj}) is strict.

\noindent
{\it 2-}
If $\alpha$ is a power of small projections, then
$k\in I\cup J\Longrightarrow\alpha _k=0$.
Thus
$$
\langle\beta ,\lambda\rangle\leq
\frac 12\sum _{k\in K}(\alpha _k -\sigma _k )
-\frac 12\sum _{k\in J}\sigma _k-\sum _{k\in I}\sigma _k
=\frac 12|\alpha |-\frac 12\delta _I^*(\sigma ).
$$
Since $\sigma$ lies in $\Sigma$, the number $\delta _I^*(\sigma )$
is nonnegative.
\end{pff}\QED

\begin{rem}
\label{equivSigma}
Assertion~{\it 1-} of Proposition~\ref{majorationsRe} is not far from being an
equivalence in the following sense.

\vskip 5pt\noindent
{\bf Claim}\hskip 5pt
\it 
For any $\alpha=(\alpha _1,\dots,\alpha _s)\in (\g R_{\geq 0})^s$,
the following are equivalent:

(1)\hskip 10pt $\forall k\in\{ 1,\dots ,s\}$,
$\Re (\lambda _k)\leq 1/2\Longrightarrow\alpha _k=0$;

(2)\hskip 10pt $\forall\beta\in (\g R_{\geq 0})^s\cap (\alpha -\Sigma )$,
$\Re\langle\beta ,\lambda\rangle <\Re\langle\alpha ,\lambda\rangle$.
\rm

\vskip 2pt
\begin{pff}
Note first that {\it (2)} is readily equivalent to
{\it (2')}: 
$\forall\sigma\in\Sigma$, $\alpha -\sigma\geq 0\Longrightarrow
\Re\langle\sigma ,\lambda\rangle> 0$.
The proof of implication {\it (1)}$\Rightarrow${\it (2')} is essentially the
same as in Proposition~\ref{majorationsRe}.
To show the contrapositive of {\it (2')}$\Rightarrow${\it (1)}, let $k$ be
such that $\Re (\lambda _k)\leq 1/2$ and $\alpha _k>0$.
Let $\sigma =\frac 12\alpha _k(2\delta _k-\delta _1)$.
It follows immediately from~(\ref{defSigma}) that $\sigma\in\Sigma$.
Furthermore, $\alpha -\sigma\geq 0$ and
$\Re\langle\sigma ,\lambda\rangle =\alpha _k(\Re(\lambda _k)-\frac 12)\leq 0$.
\end{pff}\QED

\vskip -5pt\noindent
This claim gives another ``universal'' aspect of the cone $\Sigma$ in the phase
transition phenomena between small and large processes
(see Relation (\ref{espUespQ}), Theorem~\ref{coefflemma} and  proof of
Theorem~\ref{jointmoments}).
\end{rem}

\subsubsection{Polyhedra $A_\alpha$}
\label{Aalpha-subsubsec}

Take a P\'olya process and fix a Jordan basis of linear forms
$(u_k)_{1\leq k\leq s}$.
Notations relative to this basis are defined in Section~\ref{Preliminaries}.
Remember that $\delta _k$ is the $k$-th vector of the canonical basis of
$\g R^s$.

\begin{Def}
\label{defAalpha}
For any $\alpha\in (\g Z_{\geq 0})^s$, let $A_\alpha$ be the subset of
$(\g Z_{\geq 0})^s$ defined by\begin{equation}
\label{formulaAalpha}
A_\alpha =(\alpha -D_\alpha )\cap (\g Z _{\geq 0})^s,
\end{equation}
where $D_\alpha$ is the set of $\g R_{\geq 0}$-linear combinations of all
vectors $\delta _k-\delta _{k-1}$ such that $\alpha _k\geq 1$ and
$\varepsilon _k=1$ (as in~(\ref{B-Sigma}), $\alpha -D_\alpha$ denotes
$\{\alpha -d,~d\in D_\alpha\}$).
\end{Def}

\noindent
These subsets are finite because $\beta <\beta -\delta _k+\delta _{k-1}$
for the degree-antialphabetical order defined by~(\ref{ordre}).
Geometrically speaking, $A_\alpha$ is the set of integer points of the convex
compact polyhedron of codimension $\geq 2$ defined by intersection of
$(\g R _{\geq 0})^s$ with the rational cone $\alpha -D_\alpha$ 
(and $A_\alpha$ itself is abusively called polyhedron).
When $\alpha$ is a semisimple power, all $\varepsilon _k$ vanish so that
$D_\alpha=(0)$ and $A_\alpha =\{\alpha \}$.
The polyhedra $A_\alpha$ play thus a role only for non semisimple P\'olya
processes.
Properties of $A_\alpha$'s that will be used in the sequel are listed in the
following proposition.
Definition of the differential operator $\Phi _\partial$ was given
in~(\ref{defPhipartial}).

\begin{Prop}
\label{propAalpha}
Let $\alpha\in (\g Z_{\geq 0})^s$.

1- If $\alpha$ is a semisimple power, then $A_\alpha =\{\alpha \}$.

2- If $\alpha$ is a power of large (respectively small) projections,
then every element of $A_\alpha$ is a power of large
(resp. small) projections.

3- $\alpha\in A_\alpha$ and $\Vect\{\uu ^\beta ,~\beta\in A_\alpha\}$
is $\Phi _\partial$-stable.
\end{Prop}

\begin{pff}
Justification of {\it 1-} was given just before Proposition~\ref{propAalpha}.
If $\alpha$ is a power of large (respectively small) projections,
Assertion~{\it 2-}
can be deduced from $A_\alpha$'s definition (\ref{formulaAalpha})
by induction on $\alpha$ (degree-antialphabetical order):
if $\alpha _k\geq 1$ and $\varepsilon _k=1$, then
$\alpha -\delta _k+\delta _{k-1}$ is $<\alpha$ and a power of large
(resp. small) projections.
It remains to prove {\it 3-}.
Assertion $\alpha\in A_\alpha$ is an
obvious consequence of Definition~\ref{defAalpha}.
Moreover, Jordan basis properties (see Subsection \ref{jordanbasis})
imply that $\Phi _\partial (u_1)=u_1$
and $\Phi _\partial (u_k)=\lambda _ku_k+\varepsilon _ku_{k-1}$
if $k\geq 2$, where $\varepsilon _k\in\{ 0,1\}$.
Formula~(\ref{diffualpha})
shows that for every $\beta$,
$\Phi _\partial (\uu ^\beta )-\langle\beta ,\lambda\rangle\uu ^\beta$
is linear combination of polynomials $\uu ^\gamma$, where
$\gamma =\beta -\delta _k+\delta _{k-1}$ for integers $k\geq 2$
such that $\beta _k\geq 1$ and $\varepsilon _k=1$
(hence $\lambda _k=\lambda _{k-1}$).
These considerations suffice to prove {\it 3-}.
\end{pff}\QED

\begin{Prop}
\label{propAalphalambda}
Take a P\'olya process, choose any Jordan basis and denote
$\lambda =(\lambda _1,\dots ,\lambda _s)$ its root's $s$-uple
(see~(\ref{lambda})).
Then, for every $\alpha ,\alpha '\in (\g Z_{\geq 0})^s$,
\begin{equation}
\label{equalityAalphalambda}
\alpha '\in A_\alpha\Longrightarrow
\langle\alpha ',\lambda\rangle =\langle\alpha ,\lambda\rangle
\end{equation}
\end{Prop}

\begin{pff}
If $\beta$, $\gamma$ and $k$ are like in the end of
Proposition~\ref{propAalpha}'s proof,
$\langle\gamma ,\lambda\rangle =\langle\beta ,\lambda\rangle$.
This fact leads to the result.
\end{pff}\QED

\subsection{Action of $\Phi$ on polynomials (continued)}
\label{continued-subsec}

Adopting notations of Section~\ref{cones-subsec}, we define, for any
$\alpha\in (\g Z_{\geq 0})^s$ the subspace $S'_\alpha$ of $S_\alpha$ as
\begin{equation}
\label{defS'alpha}
S'_\alpha =\Vect\{\uu ^\beta ,~\beta\in A_\alpha -\Sigma\}.
\end{equation}

\begin{Th}
\label{coefflemma}
Take a P\'olya process, choose a Jordan basis $(u_k)_{1\leq k\leq s}$ and let
$(Q_\alpha)_{\alpha\in (\g Z_{\geq 0})^s}$ be the corresponding reduced
polynomials.

\vskip 3pt
1- For any $\alpha\in (\g Z_{\geq 0})^s$,
\begin{equation}
\label{S'alphastable}
S'_\alpha {\rm ~is~\Phi-stable};
\end{equation}
\begin{equation}
\label{span=span}
S'_\alpha =\Vect \{Q_\beta ,~\beta\in A_\alpha -\Sigma\}.
\end{equation}

2- If $\alpha$ is a power of large projections, then
$\Phi (Q_\alpha )\in\Vect \{Q_\beta ,~\beta\in A_\alpha\}$.
\end{Th}

\noindent
Consequently, Relations (\ref{qalphabeta})
can be refined in the general case by means of~(\ref{span=span})
and Proposition~\ref{proprietesQalpha}; this provides straightforwardly
\begin{equation}
\label{refined}
\uu ^\alpha =Q_\alpha
+\sum _{\beta\in A_\alpha -\Sigma ,
~\langle\beta ,\lambda\rangle\neq\langle \alpha ,\lambda\rangle}
q_{\alpha ,\beta}Q_\beta .
\end{equation}

\begin{pff}
{\it 1-}
We first prove that $S'_\alpha$ is $\Phi$-stable.
Let $\beta\in A_\alpha -\Sigma$;
let $\alpha '\in A_\alpha$ and $\sigma\in\Sigma$ such that
$\beta =\alpha '-\sigma$.
We show that both $\Phi _\partial (\uu ^\beta )$ and
$(\Phi -\Phi _\partial )(\uu ^\beta )$ belong to $S'_\alpha$.

\noindent
$\bullet$
As in Proposition~\ref{propAalpha}'s proof,
$\Phi _\partial (\uu ^\beta )-\langle\beta ,\lambda\rangle\uu ^\beta$
is a linear combination of polynomials $\uu ^\gamma$, where
$\gamma =\beta -\delta _k+\delta _{k-1}$ for integers $k\geq 2$
such that $\beta _k\geq 1$ and $\varepsilon _k=1$
(hence $\lambda _k=\lambda _{k-1}$).
If $\gamma$ is such a power, we claim that
$\gamma\in A_\alpha -\Sigma$, which shows that 
$\Phi _\partial (\uu ^\beta )\in S'_\alpha$.
If $\alpha '_k\geq 1$, just write $\gamma =\alpha ''-\sigma$
where $\alpha ''=\alpha '-\delta _k+\delta _{k-1}\in A_\alpha$.
If $\alpha '_k=0$, this $\alpha ''$ is not in $A_\alpha$ because it is
not nonnegative;
in this case, write $\gamma =\alpha '-\sigma '$ where
$\sigma '=\sigma +\delta _k-\delta _{k-1}$.
It suffices to show that $\sigma '\in\Sigma$.
Let $I$ be a proper subset of $\{ 1,\dots ,s\}$, that gives the equation
of a face of $\Sigma$ (see Corollary~\ref{facesSigma}).
If $k\in I$ or $k-1\notin I$, then
$\delta _I^*(\sigma ')\geq 0$ because $\delta _I^*(\sigma )\geq 0$.
If $k\notin I$ and $k-1\in I$, then
$\delta _I^*(\sigma ')=\delta _I^*(\sigma )-1
=\delta _{I\cup \{k\} }^*(\sigma )-\sigma _k-1$;
but $\sigma _k=-\beta _k\leq -1$ because $\beta _k\geq 1$.
Thus $\delta _I^*(\sigma ')\geq\delta _{I\cup \{k\} }^*(\sigma )\geq 0$
since $\sigma\in\Sigma$
(note that this last inequality is true even if $I\cup\{ k\} =\{ 1,\dots ,s\}$
because $(1,\dots ,1)\in\check\Sigma$, see Remark~\ref{deltavide}).

\noindent
$\bullet$
Taylor formula implies that
$(\Phi -\Phi _\partial)(\uu ^\beta )$ is linear combination of polynomials
$\uu ^\gamma =\uu ^{\beta -\eta +\delta _k}$ with $1\leq k\leq s$,
$\eta \geq 0$, $|\eta |\geq 2$, $\beta -\eta\geq 0$
(the $\eta$-terms correspond to partial derivatives of order $\geq 2$ of
$\uu ^\beta$, the $\delta _k$-terms come from the expansion of linear
forms $l_k$ in the Jordan basis $(u_k)_k$).
If $\gamma =\beta -\eta +\delta _k$ is such a power and if
$\delta _I^*$
is the equation of any one of the defining hyperplanes of $\Sigma$
where $I$ is a proper subset of $\{1,\dots ,s\}$
(see Corollary \ref{facesSigma}), then
$$
\delta _I^*(\beta -\gamma)=\delta _I^*(\eta -\delta _k)
=|\eta|-1+\sum _{i\in I}\eta _i-\11 _{k\in I}
\geq 1-\11 _{k\in I}+\sum _{i\in I}\eta _i\geq 0.
$$
This proves that $(\Phi -\Phi _\partial)(\uu ^\beta )$
belongs to
$\Vect\{\uu ^\gamma,~|\gamma |\leq |\beta |-1,~\gamma \in\beta -\Sigma\}$.
If $\gamma =\beta -\sigma '\in\beta -\Sigma$, then
$\gamma=\alpha '-(\sigma +\sigma ')\in\alpha '-\Sigma$
because the cone $\Sigma$ is stable under addition.
This shows that $(\Phi -\Phi _\partial)(\uu ^\beta )\in S'_\alpha$
(see Figures \ref{Fi:sigma2} and \ref{Fi:sigma3} for a representation of
powers $\beta -\gamma =\eta -\delta _k$ that appear in this computation).

Thus, $S'_\alpha$
is a $\Phi$-stable subspace of $S_\alpha$.
For any $\beta\in A_\alpha -\Sigma$, the projection of $\uu ^\beta$ on
$S'_\alpha\cap\ker (\Phi -\langle\beta ,\lambda\rangle )^\infty$
parallel to
$S'_\alpha\cap\bigoplus _{z\neq\langle\beta ,\lambda\rangle}
\ker (\Phi -z)^\infty$
equals $Q_\beta$ because of the unicity of the decomposition on
characteristic spaces.
Hence $Q_\beta\in S'_\alpha$
(another way to show that fact consists in noting that the projections on these
characteristic spaces are polynomial functions of the restriction of 
$\Phi$ to $S'_\alpha$).
Thus, $\Vect \{Q_\beta ,~\beta\in A_\alpha -\Sigma\}$ is a subspace of
$S'_\alpha$;
as these two subspaces have the same finite dimension, they are equal.
The proof of {\it 1-} is complete.

{\it 2-}
Because of Assertion {\it (5)} in Proposition~\ref{proprietesQalpha} and of
the $\Phi$-stability of $S'_\alpha$,
$$
\Phi (Q_\alpha )\in
\Vect \{Q_\beta ,~\alpha <\beta ,~\beta\in A_\alpha -\Sigma ,
~\langle\beta ,\lambda\rangle =\langle\alpha ,\lambda\rangle\}.
$$
Conclude with {\it 1-} of Proposition~\ref{majorationsRe} and
Proposition~\ref{propAalphalambda}.
\end{pff}\QED

\begin{rem}
{\bf (On the minimality of cone $\Sigma$ and polyhedra $A_\alpha$)}
\label{remarkSigma}

Cone $\Sigma$ appears in a natural way
in the proof of Theorem \ref{coefflemma} to ensure the $\Phi$-stability
of a (minimal) subspace that contains some given $\uu ^\alpha$.
Indeed, suppose for simplicity that $\alpha$ is a semisimple power.
Then $\Phi (\uu ^\alpha )$ is the sum of
$\langle\alpha ,\lambda\rangle\uu ^\alpha$ and of a linear combination
of polynomials $\uu ^{\alpha -\eta +\delta _k}$ where $\eta\geq 0$,
$|\eta |\geq 2$ and $1\leq k\leq s$.
The iterations of $\Phi$ on such polynomials force to consider the least
(for inclusion) set of powers that contains these $\eta -\delta _k$
and that is stable under addition (and contains zero);
this least set is $\Sigma$.
For an illustration of this fact, see Figures~\ref{Fi:sigma2}
and~\ref{Fi:sigma3}.
If $\alpha$ is not semisimple, the situation is complicated by other
powers $\alpha '$ of same total degree and leads to consider the polyhedron
$A_\alpha$.

Suppose that $B\subseteq (\g Z_{\geq 0})^s$ is such that $\alpha\in B$
and $\Vect\{\uu ^\beta ,~\beta\in B\}$ is $\Phi _\partial$-stable.
As in the end of Proposition~\ref{propAalpha}'s proof,
Formula~(\ref{diffualpha}) shows that $A_\alpha\subseteq B$, so that
$A_\alpha$ is minimal for properties {\it 3-} of Proposition~\ref{propAalpha}.
These properties are necessary to imply Formulae~(\ref{span=span})
and ~(\ref{refined}).
The optimality of $A_\alpha$, announced in the introduction of
Subsection~\ref{cones-subsec}, consists in that fact.

A shorter proof of Theorem~\ref{jointmoments} can nevertheless be given without
considering~$\Sigma$.
See Remark~\ref{GammaReferee}.
\end{rem}

\masection{Proof of Theorem~\ref{jointmoments}, asymptotics of moments}
\label{jointmoments-sec}

The proof of Theorem~\ref{jointmoments} relies on Formula~(\ref{refined}) in
which the expectation of the value at $X_n$ is taken, providing
\begin{equation}
\label{finalexpansion}
E\uu ^\alpha (X_n)=EQ_\alpha (X_n)
+\sum _{\beta\in A_\alpha -\Sigma ,
~\langle\beta ,\lambda\rangle\neq\langle \alpha ,\lambda\rangle}
q_{\alpha ,\beta}EQ_\beta (X_n).
\end{equation}
We first use Corollary~\ref{fjordan} to make asymptotics of reduced moments
$EQ_\beta (X_n)$ precise in Proposition~\ref{reducedmoments}
before giving a proof of Theorem~\ref{jointmoments}.
Numbers $\nu _\alpha$ were defined by~(\ref{defnualpha}) and appear naturally
in Theorem~\ref{jointmoments};
Subsection~\ref{nualpha-subsec} is devoted to their computation in the case
of powers of large projections.

\subsection{Proof of Theorem~\ref{jointmoments}}
\label{jointmoments-subec}

\begin{Prop}
\label{reducedmoments}
{\bf (Asymptotics of reduced moments)}
Let $\alpha\in (\g Z _{\geq 0})^s$.

1- If $\nu _\alpha =0$, i.e. if $Q_\alpha$ is eigenfunction of $\Phi$, then,
as $n$ tends to infinity,
$$
EQ_\alpha (X_n)=n^{\langle\alpha ,\lambda\rangle}
\frac {\Gamma (\tau _1)}{\Gamma (\tau _1+\langle\alpha ,\lambda\rangle )}
Q_\alpha (X_1)
+O(n^{\langle\alpha ,\lambda\rangle -1}).
$$

2- If $\nu _\alpha\geq 1$, then, as $n$ tends to infinity,
$$\begin{array}{rc}
EQ_\alpha (X_n)&\hskip -5pt =\hskip 3pt\displaystyle
\frac {n^{\langle\alpha ,\lambda\rangle}\log ^{\nu _\alpha}n}{\nu _\alpha !}
\frac{\Gamma (\tau _1)}{\Gamma (\tau _1+\langle\alpha ,\lambda\rangle )}
(\Phi -\langle\alpha ,\lambda\rangle )^{\nu _\alpha}(Q_ \alpha )(X_1)\\
&\\
&\hskip 80pt +~O(n^{\langle\alpha ,\lambda\rangle}\log ^{\nu _\alpha -1}n).\\
\end{array}$$
\end{Prop}

\vskip 5pt
\begin{pff}
$Q_\alpha$ belongs to the $\Phi$-stable subspace
${\cal S}(\alpha )=\ker (\Phi -\langle\alpha ,\lambda\rangle )^\infty$
and the operator induced by $\Phi$ on ${\cal S}(\alpha )$ is the sum of
$\langle\alpha ,\lambda\rangle\Id _{{\cal S}(\alpha )}$ and of the nilpotent
operator induced on ${\cal S}(\alpha )$ by
$\Phi -\langle\alpha ,\lambda\rangle$.
This facts being considered, Proposition~\ref{reducedmoments} is a
straightforward consequence of Corollary~\ref{fjordan}.
\end{pff}\QED

\begin{rem}
Even if $\alpha$ is a semisimple power, $Q_\alpha$ may not be an eigenfunction
of $\Phi$.
This happens only if
$\langle\alpha ,\lambda\rangle =\langle\beta ,\lambda\rangle$
for some $\beta <\alpha$ because this implies that both $Q_\alpha$ and
$Q_\beta$ are in the same characteristic subspace of $\Phi$.
When all roots $\lambda _k$ are incommensurable, {\it i.e.} whenever they admit
no non trivial linear relation with rational coefficients, all numbers
$\langle\alpha ,\lambda\rangle$ are distinct, so that every $Q_\alpha$ is
an eigenfunction of $\Phi$.

\end{rem}

\noindent
{\sc Proof of Theorem~\ref{jointmoments}.}\ 
Proposition \ref{reducedmoments} asserts in particular that for all
$\beta \leq \alpha$, there exists $\nu\geq 0$ such that
$EQ_\beta (X_n)\in
O\left( n^{\Re\langle\beta ,\lambda\rangle}\log ^\nu n\right)$.

\noindent
{\it 1-}
If $\alpha$ is a power of small projections, then any $\alpha '\in A_\alpha$
satisfies $|\alpha '|=|\alpha |$ and is a power of small projections
(Proposition~\ref{propAalpha}).
Hence $n^{\langle\beta ,\lambda\rangle}\in O\left( n^{|\alpha |/2}\right)$
if $\beta\in\alpha '-\Sigma$, as can be deduced from Proposition
\ref{majorationsRe}.

\noindent
{\it 2-}
If $\alpha$ is a power of large projections, then
any $\alpha '\in A_\alpha$ satisfies
$\langle\alpha ',\lambda\rangle =\langle\alpha ,\lambda\rangle$
and  is power of large projections.
Hence, for every $\beta\in \alpha '-\Sigma$,
$\alpha '=\beta$ or
$\Re\langle\beta ,\lambda\rangle <\Re\langle\alpha ,\lambda\rangle$
(Proposition \ref{majorationsRe} and Proposition \ref{propAalphalambda}).
Thus, Formula~(\ref{finalexpansion}), implies that
$E\uu ^\alpha (X_n)=EQ_{\alpha }(X_n)
+o\left( n^{\Re\langle\alpha ,\lambda\rangle }\right)$
has the required asymptotics as it can be deduced from
Proposition \ref{reducedmoments}.

\noindent
{\it 3-}
Moreover, if $\alpha$ is a semisimple power of large projections, then
$\nu _\alpha =0$
(consequence of Theorem~\ref{coefflemma}~({\it 2-}),
Proposition~\ref{propAalpha}~({\it 1-}) and $\nu _\alpha$'s
definition~(\ref{defnualpha})).
We conclude with Proposition \ref{reducedmoments}.
\QED

\begin{rem}
\label{constantc}
In Theorem~\ref{jointmoments}, Assertion~{\it 3-} is a particular case of
Assertion~{\it 2-},
the constant named $c$ appearing in {\it 2-} being
$$
c=\frac 1{\nu _\alpha !}
\frac{\Gamma (\tau _1)}{\Gamma (\tau _1+\langle\alpha ,\lambda\rangle )}
(\Phi -\langle\alpha ,\lambda\rangle )^{\nu _\alpha}(Q_ \alpha )(X_1).
$$
This follows from Proposition \ref{reducedmoments} and from
Theorem~\ref{jointmoments}'s proof.
\end{rem}

\begin{rem}
More precision on the small $o$ of Assertion~{\it 3-}
in Theorem~\ref{jointmoments} can be deduced from its proof:
one can replace it by $O\left( n^a\right)$ where
$$
a=\max\big(\{\Re\langle \beta ,\lambda \rangle,
~\beta\neq\alpha ,~\beta\in\alpha -\Sigma\}
\cup 
\{\Re\langle \alpha ,\lambda \rangle -1\}\big) .
$$
\end{rem}

\begin{rem}
\label{GammaReferee}
One can give a shorter proof of Theorem~\ref{jointmoments}
without explicitely considering $\Sigma$.
The following arguments, that provide a self-sufficient independent proof,
 have been suggested by the anonymous referee.
For any $\beta =(\beta _1,\dots ,\beta _s)\in (\g Z_{\geq 0})^s$ and any
$z\in\g C$, let's define
$\beta _{[z]}=\sum _{1\leq k\leq s,~\lambda _k=z}\beta _k$
(this number being $0$ by convention if there is no $k$ such that
$\lambda _k=z$), and let
$$
\Gamma =\left\{
\beta\in (\g Z_{\geq 0})^s,~
\sum _{z,~\beta _{[z]}\geq 0}\beta _{[z]}
+2\sum _{z,~\beta _{[z]}\leq 0}\beta _{[z]}\geq 0
\right\}.
$$
Then
$\Gamma$ is stable under addition (it is a convex cone) and the subspace
$F_\alpha=
\Vect\{\uu ^\beta,~\beta\in (\g Z_{\geq 0})^s\cap (\alpha -\Gamma )\}$
is thus $\Phi$-stable as can be seen from
Equation~(\ref{diffualpha}) and from the proof of Proposition~\ref{phiordre},
page \pageref{preuvephiordre}.
The fact that $F_\alpha$ is $\Phi$-stable implies that $Q_\beta\in F_\alpha$
for all $\beta\in (\g Z_+)^s\cap (\alpha -\Gamma)$;
indeed, $Q_\beta$ can be seen as the projection of $\uu ^\beta$ on
$F_\alpha\cap\ker (\Phi -\langle\beta ,\lambda\rangle )^\infty$ parallel to
$F_\alpha\cap
\bigoplus _{z\neq\langle\beta ,\lambda\rangle}\ker (\Phi -z)^\infty$,
because of the unicity of the decomposition on characteristic spaces.
Therefore,
$F_\alpha=
\Vect\{ Q_\beta,~\beta\in (\g Z_{\geq 0})^s\cap (\alpha -\Gamma )\}$.
Besides, $\Gamma$ has the following two properties:

{\it (a)}
if $\alpha$ is a power of large projections, then $\beta\in\alpha -\Gamma$
and $\langle\beta ,\lambda\rangle\neq\langle\alpha ,\lambda\rangle$
imply $\Re\langle\beta ,\lambda\rangle <\Re\langle\alpha ,\lambda\rangle$;

{\it (b)}
if $\alpha$ is a power of small projections, then
$\Re\langle\beta ,\lambda\rangle\leq |\alpha |/2$
for any $\beta\in (\g Z_{\geq 0})^s\cap (\alpha -\Gamma )$.

\noindent
[Proof of {\it (a)}:
if $\gamma =\alpha -\beta$ is such that $(\gamma _{[z]})_{z\in\g C}\neq 0$,
then
$$
\begin{array}{rcl}
\Re\langle\gamma ,\lambda\rangle &
=&\displaystyle
\sum _{\gamma _k\geq 0}\gamma _k\Re\lambda _k
+\sum _{\gamma _k\leq 0}\gamma _k\Re\lambda _k\\
&>&\displaystyle
\frac 12\sum _{\gamma _k\geq 0}\gamma _k+\sum _{\gamma _k\leq 0}\gamma _k
\geq 0.
\end{array}
$$
Assertion {\it (b)} follows from a similar argument.]

\noindent
These two properties suffice to show Assertions~{\it 1-} and~{\it 2-} in
Theorem~\ref{jointmoments}.
Likewise, Theorem~\ref{coefflemma}~{\it 2-} and
Theorem~\ref{jointmoments}~{\it 3-} can be obtained by the following
refinement of the argument:
$\Vect\{\uu ^\beta ,~\beta\in A_\alpha ~{\rm or~}\beta =\alpha -\sigma,~
\sigma\in\Gamma,~(\sigma _{[z]})_{z\in\g C}\neq 0\}$
is $\Phi$-stable.

In any case, $\Sigma\subseteq\Gamma$.
When all roots of the process are distincts, then $\Gamma =\Sigma$.
Otherwise, the cone $\Gamma$ (more precisely $\Gamma\otimes\g R$) is not
strictly convex (it contains a nonzero vector subspace of $\g R^s$).
Consequently, since $\Sigma$ is strictly convex, $\Gamma =\Sigma$
if and only if all roots of the process are distinct.
Compare this fact with the minimality of $\Sigma$ asserted in
Remark~\ref{remarkSigma}.
\end{rem}

\subsection{Computation of nilpotence indices $\nu _\alpha$}
\label{nualpha-subsec}

For any $\alpha\in (\g Z_{\geq 0})^s$, the number $\nu _\alpha$ has been
defined in Section~\ref{results-sec} as the nilpotence index of $Q_\alpha$
for $\Phi$, {\it i.e.}
$$
\nu _\alpha=\max\{ p\geq 0,~(\Phi -\langle\alpha ,\lambda\rangle )^p(Q_\alpha )
\neq 0\}.
$$
It appears in the expression of the leading term of the $\uu$-moment
$E\uu ^\alpha (X_n)$ as $n$ tends to infinity, when $\alpha$ is a power
of large projections (Theorem~\ref{jointmoments}).

In several problems where these moments' asymptotics are needed, it is useful
to compute them explicitely.
To this end, iterations of the finite difference operator $\Phi$ are not
easily handled;
furthermore, a calculation of $\nu _\alpha$ from its definition supposes that
the reduced polynomial $Q_\alpha$ has already been computed.
These two facts make a direct computation of $\nu _\alpha$ rather intricate.
Whenever $\alpha$ is a power of large projections, 
Proposition~\ref{nualpha} asserts that $\nu _\alpha$ is the nilpotence index of
$\uu ^\alpha$ for the differential operator $\Phi _\partial$,
making its computation much easier.

\begin{Prop}
\label{nualpha}
{\bf (Computation of $\nu _\alpha$)}

1- If $\alpha$ is a power of large projections, then
$$
\nu _\alpha =\max\{ q\geq 0,
~(\Phi _\partial -\langle\alpha ,\lambda\rangle )^q(\uu ^\alpha )\neq 0\} .
$$

2- If $\alpha =(\alpha _1,\dots ,\alpha _s)$ is a monogenic power of large
projections whose support is contained in the monogenic
block of indices $J=\{m,\dots ,m+r\}$ ($r\geq 0$), then
$$
\nu _\alpha =\sum _{k=0}^rk\alpha _{m+k}.
$$
\end{Prop}

\begin{pff}
{\it 1-} There is nothing to prove if $\alpha$ is a semisimple power.
We thus suppose that $\alpha$ is a large and not semisimple power.
We denote by $k_\alpha$ the index
$k_\alpha =\min\{k\geq 3, ~\alpha _k\geq 1,~\varepsilon _k=1\}$
and $p(\alpha )$ the element of $A_\alpha$ defined by
$p(\alpha )=\alpha -\delta _{k_\alpha}+\delta _{k_\alpha -1}$;
it is the predecessor of $\alpha$ for the degree-antialphabetical
order restricted to $A_\alpha$.
As a direct computation of $\Phi _\partial(\uu ^\alpha )$ shows
(see (\ref{diffualpha}) in the proof of Proposition \ref{phidiffpropre}),
Proposition~\ref{propAalpha} implies that
\begin{equation}
\label{predPhipartialu}
(\Phi _\partial -\langle\alpha ,\lambda\rangle )(\uu ^\alpha )
-\alpha _{k_\alpha}\uu ^{p(\alpha )}
\in\Vect\{\uu ^\beta ,~\beta\in A_\alpha ,~\beta <p(\alpha )\} .
\end{equation}
We claim that
$
(\Phi -\langle\alpha ,\lambda\rangle )(Q_\alpha )
-\alpha _{k_\alpha}Q_{p(\alpha )}
\in\Vect\{Q_\beta ,~\beta <p(\alpha )\}
$
(proof just below).
With Theorem~\ref{coefflemma},
as $\alpha$ is a power of large projections, this implies that
\begin{equation}
\label{predPhiQ}
(\Phi -\langle\alpha ,\lambda\rangle )(Q_\alpha )
-\alpha _{k_\alpha}Q_{p(\alpha )}
\in\Vect\{Q_\beta ,~\beta\in A_\alpha ,~\beta <p(\alpha )\} .
\end{equation}
Assertions~(\ref{predPhipartialu}) and~(\ref{predPhiQ}) are then sufficient
to show that
$$
\max\{ q\geq 0,~(\Phi -\langle\alpha ,\lambda\rangle )^q(Q_\alpha )\neq 0\}
=
\max\{ q\geq 0,
~(\Phi _\partial -\langle\alpha ,\lambda\rangle )^q(\uu ^\alpha )\neq 0\},
$$
this number being equal to
$
\min\{ q\geq 0,~p^{[q]}(\alpha )~is~a~semisimple~power\}
$
(notation $p^{[q]}$ denotes the composition $p\circ\dots p$ iterated
$q$ times and $p^{[0]}(\alpha )=\alpha$).

It remains to prove that
$(\Phi -\langle\alpha ,\lambda\rangle )(Q_\alpha )
-\alpha _{k_\alpha}Q_{p(\alpha )}
\in\Vect\{Q_\beta ,~\beta <p(\alpha )\}$.
Note first that $|p(\alpha )|=|\alpha |$, so that $\beta <p(\alpha )$ as soon
as $|\beta |\leq |\alpha |-1$.
Since $\Phi -\Phi _\partial$ let the degree fall down (Taylor formula),
\begin{equation}
\label{tructruc}
(\Phi -\langle\alpha ,\lambda\rangle )(Q_\alpha )
\in (\Phi _\partial -\langle\alpha ,\lambda\rangle )(Q_\alpha )
+\Vect\{\uu ^\beta ,~|\beta |\leq |\alpha |-1\}
.
\end{equation}
The only point of $\Sigma$ having a nonpositive degree ($|.|$) is zero, so that
Theorem~\ref{coefflemma} leads to
$$
Q_\alpha \in\uu ^\alpha +\Vect\{\uu ^\beta ,~|\beta |\leq |\alpha |-1\}
+\Vect\{\uu ^\beta ,~\beta <\alpha ,~\beta\in A_\alpha\}.
$$
Taking the image by $\Phi _\partial -\langle\alpha ,\lambda\rangle$
of this last relation leads, using~(\ref{tructruc}) to
$$
(\Phi -\langle\alpha ,\lambda\rangle )(Q_\alpha )
\in (\Phi _\partial -\langle\alpha ,\lambda\rangle )(\uu ^\alpha )
+\Vect\{\uu ^\beta ,~\beta <p(\alpha )\}
.$$
Because of Assertion (\ref{predPhipartialu}),
$$
(\Phi -\langle\alpha ,\lambda\rangle )(Q_\alpha )
\in\alpha _{k_\alpha}\uu ^{p(\alpha )}
+\Vect\{\uu ^\beta ,~\beta <p(\alpha )\}
.$$
The conclusion follows from Proposition~\ref{proprietesQalpha}
(Assertions {\it 2-} and {\it 4-}).

{\it 2-}
The total degree $|\alpha |$ being fixed, we proceed by induction on $\alpha$.
If $\alpha$ is semisimple, $\alpha =|\alpha |\delta _m$ and
$\nu _\alpha =0$; there is nothing to prove.
If $\alpha$ is not semisimple, the computation of
$(\Phi _\partial -\langle\alpha ,\lambda\rangle )(\uu ^\alpha )$
shows that
$$
\nu _\alpha =1+\max\{\nu _{\alpha -\delta _k+\delta _{k-1}},
~m+1\leq k\leq m+r,~\alpha _k\geq 1\} .
$$
All these $\alpha -\delta _k+\delta _{k-1}$ are $<\alpha$ and have total degree
$|\alpha |$;
by induction, they all have the same $\nu$,
this number being $-1+\sum _{0\leq k\leq r}k\alpha _{m+k}$.
The formula for $\nu _\alpha$ is proven.
\end{pff}\QED

\masection{Proofs of Theorems~\ref{lpss} and~\ref{lpnss},
asymptotics of large processes}
\label{main-sec}

\noindent
{\sc Proof of Theorem~\ref{lpss}.}\ 
We adopt the notations of Section \ref{Preliminaries}.
Let's denote $\pi =\sum_{2\leq k\leq r}\pi _k$
and $\pi '=\sum_{k\geq r+1}\pi _k$;
the random vector $X_n$ splits into the sum
\begin{equation}
\label{trisection}
X_n=\pi _1X_n+Y_n+Z_n,
\end{equation}
where $Y_n=\pi X_n$ and $Z_n=\pi 'X_n$.

\noindent
$\bullet$
First term $\pi _1X_n$

Definitions of the Jordan basis $(u_k)_k$ of linear forms and of its dual basis
$(v_k)_k$ of vectors imply readily that $\pi _1(v)=u_1(v)v_1$ for any vector
$v$.
Thus, because of Relation (\ref{balance}),
$\pi _1X_n=nv_1+O(1)$ as $n$ tends to infinity;
this projection is non random.

\noindent
$\bullet$
Second term $Y_n$

As follows from (\ref{projections}), $Y_n=\sum _{k=2}^ru_k(X_n)v_k$.
Take any $k\in\{2,\dots ,r\}$.
Computation of the expectation of $u_k(X_{n+1})$ conditionally
to the state at time~$n$ gives
$E^{{\cal F}_n}u_k(X_{n+1})=(1+\lambda _k/(n+\tau _1-1))u_k(X_n)$
for any positive integer~$n$
(see~(\ref{espeCondi}), $u_k$ is an eigenform of the process);
this implies that
$(\gamma _{\tau _1,n}(\lambda _k)^{-1}u_k(X_n))_n$
is a martingale
(one can divide by $\gamma _{\tau _1,n}(\lambda _k)$ because $\lambda _k$
is not a negative integer).
As $\overline{u_k}$ (complex conjugacy) is an eigenform associated to the
eigenvalue $\overline{\lambda _k}$, it is a linear combination of eigenforms
$u_l$'s, all associated with $\overline{\lambda _k}$.
Thus, if $q\geq 1$ is an integer, $|u_k^{2q}|$ is a linear combination of
polynomials $\uu ^\alpha$'s for some suitable semisimple powers of large
projections $\alpha$'s such that $\langle\alpha ,\lambda\rangle =2q\sigma _2$.
This implies, thanks to Theorem \ref{jointmoments}, that
$$
E|u_k^{2q}(X_n)|=O(n^{2q\sigma _2}).
$$
Note that this is valid even if $\lambda _k$ is real.
The martingales
$\gamma _{\tau _1,n}(\lambda _k)^{-1}u_k(X_n)$ are consequently all convergent
in every ${\rm L}^p$ space, $p\geq 1$.

For every $k\in\{2,\dots ,r\}$, let $W_k$ be the (complex) random variable
defined by
$$
W_k=\lim _{n\to+\infty}
\frac{u_k(X_n)}{\gamma _{\tau _1,n}(\lambda _k)}
\frac{\Gamma (\tau _1)}{\Gamma (\tau _1+\lambda _k)}
=\lim _{n\to+\infty}u_k\left( X_n/n^{\lambda _k}\right)
$$
the second equality coming from Stirling's asymptotics
as $n$ tends to infinity:
$$
\gamma _{\tau _1,n}(\lambda)\Gamma (\tau _1+\lambda)
=\Gamma (\tau _1)n^\lambda (1+o(1)),
$$
for every $\lambda\notin -\tau _1+\g Z_{\leq 0}$.
This shows that
$Y_n=\sum _{2\leq k\leq r}n^{\lambda _k}W_kv_k+o(n^{\sigma _2})$,
the small $o$ being almost sure and in ${\rm L}^p$ for
any $p\geq 1$.

Computation of joint moments' limits:
if $\alpha =(0,\alpha _2,\dots ,\alpha _r,0,\dots )\in(\g Z_{\geq 0})^s$,
$\alpha$ is a semisimple power of large projections
and if one denotes $W^\alpha =\prod _{k}W_k^{\alpha _k}$, Theorem
\ref{jointmoments} implies that
$$
EW^\alpha =\lim _{n\to +\infty}
\frac1{n^{\langle\alpha,\lambda\rangle}}E\uu ^\alpha (X_n)
=\frac{\Gamma (\tau _1)}{\Gamma (\tau _1+\langle\alpha ,\lambda\rangle )}
Q_\alpha (X_1).
$$

\noindent
$\bullet$ Third term $Z_n$

$Z_n=\sum _{k\geq r+1}u_k(X_n)v_k$.
We show that $n^{-\sigma _2}u_k(X_n)$ converges to zero almost surely and
in every ${\rm L}^p$ space ($p\geq 1$), for every $k\geq r+1$.
Take any $k\geq r+1$ and any integer $q\geq 1$.
As above, $\overline{u_k}$ is a linear combination of $u_l$'s, all associated
with the root $\overline{\lambda _k}$ (even if $u_k$ and the $u_l$'s are not
necessarily eigenforms).
Thus, $|u_k^{2q}|$ is a linear combination of polynomials $\uu ^\alpha$'s for
some suitable $\alpha$'s that are powers of large (respectively small)
projections if $\Re\lambda _k>1/2$ (\emph{resp.} if $\Re\lambda _k\leq 1/2$),
such that $\langle\alpha ,\lambda\rangle =2q\Re\lambda _k$.
Because of Theorem \ref{jointmoments}, this implies in any case that
$E|u_k^{2p}(X_n)|\in o(n^{2p\sigma _2})$, which gives the ${\rm L}^p$
convergence.
Furthermore, let $p$ be any positive integer such that
$1/p< 2(\sigma _2-\Re\lambda _k)$ if $\Re\lambda _k>1/2$
or such that $1/p< 2\sigma _2-1$ if not;
for such a $p$, the series
$$
\sum _nE\left| \frac1{n^{\sigma _2}}u_k(X_n)\right| ^{2p}
$$
converges.
The almost sure convergence to zero of $n^{-\sigma _2}u_k(X_n)$ follows
thus from the almost sure convergence of the series of nonnegative random
variables
$$
\sum _n\left| \frac1{n^{\sigma _2}}u_k(X_n)\right| ^{2p},
$$
and the proof of Theorem~\ref{lpss} is complete.
\QED

\vskip 10pt
\noindent
{\sc Proof of Theorem~\ref{lpnss}.}\ 
We adopt the notations of Section \ref{Preliminaries}.
For any monogenic block of indices $J$, we denote by $\pi _J$
the projection $\pi _J=\sum _{k\in J}\pi _k$.

\noindent
$\bullet$
{\bf Claim}\hskip 5pt
{\it If $J$ is a monogenic block of indices associated with a root $\lambda$
having a real part $\sigma >1/2$, then $\gamma _{\tau _1,n}(\pi _JA)$
is invertible and
$\gamma _{\tau _1,n}(\pi _JA)^{-1}\pi _JX_n$ is a martingale that
converges in ${\rm L}^p$ for every $p\geq 1$ (thus almost surely).
If $M_J$ denotes the limit of this martingale and if $\nu=\# J -1$, then
\begin{equation}
\label{asymptpiJXn}
\pi _JX_n=
\frac {n^\lambda\log ^\nu n}{\nu !}
\frac {\Gamma (\tau _1)}{\Gamma (\tau _1+\lambda )}
u_{\min J}(M_J)v_{\max J}+o\left( n^\sigma\log ^\nu n\right)
\end{equation}
as $n$ tends to infinity, the small $o$ being almost sure and
in ${\rm L}^p$ for every $p\geq 1$.
Furthermore, almost surely and in ${\rm L}^p$ for every $p\geq 1$,
\begin{equation}
\label{u2MJ}
u_{\min J}(M_J)
=\frac{\Gamma (\tau _1+\lambda )}{\Gamma (\tau _1)}
\times
\lim _{n\to\infty}\frac{u_{\min J}(X_n)}{n^\lambda }.
\end{equation}
}

\vskip 3pt\noindent
$\bullet$
Proof of the claim.
The endomorphism
$\gamma _{\tau _1,n}(\pi _JA)$ is invertible because every
$\Id +\pi _JA/(k+\tau _1-1)$ is
(its unique eigenvalue has a real part $>1$).
Since $J$ is a monogenic block of indices, $A$ and $\pi _J$ commute.
Thus $M_n=\gamma _{\tau _1,n}(\pi _JA)^{-1}\pi _JX_n$ is a martingale
(see (\ref{espeCondi})
with $f=\pi _J$ and Remark~\ref{flinear}
in Subsection \ref{phimoments-subsec}).
We show that for any $k\in J$, the quadratic variation of the martingale
$u_k(M_n)$ is almost surely bounded, which suffices, thanks to
Burkholder's Inequality for discrete time martingales
(see \cite{HallHeyde} for example),
to ensure that the
projection $u_k(M_n)$ is bounded in ${\rm L}^p$ for every $p\geq 1$, hence
the validity of the convergence part of the claim.

Without loss of generality, we can assume for simplicity that
$J=\{2,\dots ,\nu +2\}$.
If one denotes $N=\pi _J(A-\lambda)$, then $N$ commutes with $A$ and satisfies
$N^\nu \neq 0$ and $N^{\nu +1}=0$;
furthermore, elementary considerations on $A$, the $u_k$'s and the $v_k$'s
show that for any nonnegative integer $q$ and for any
$k\in\{2,\dots ,\nu +2\}$, one has
$N^q\pi _k=u_kv_{k+q}$ if $k+q\leq \nu +2$ and $N^q\pi _k=0$ is
$k+q\geq \nu +3$.
In particular, for any $q$, one can write
$N^q=N^q(\sum _{k\in J}\pi _k)=\sum _{q+2\leq k\leq \nu +2}u_{k-q}v_k$
(with the convention $N^0=\pi _J$).
Hence, if $\beta _n=1/\gamma _{\tau _1,n}$
(as formal series or rational fraction; we omit the parameter $\tau _1$ for
simplicity of notation),
Taylor formula leads to
\begin{equation}
\label{gamman-1}
\begin{array}{rcl}
M_n
&=&\displaystyle\beta _n(\lambda +N)\pi _JX_n
=\sum _{q=0}^\nu \frac 1{q!}\beta _n^{(q)}(\lambda )N ^qX_n\\
&=&\displaystyle\sum _{k=2}^{\nu +2}
\left(\sum _{q=0}^{k-2}\frac 1{q!}\beta _n^{(q)}(\lambda )u_{k-q}(X_n)\right)
v_k.
\end{array}
\end{equation}
Thus, for any $k\in\{2,\dots ,\nu +2\}$, one has
$u_k(M_n)=\sum _{q=0}^{k-2}\frac 1{q!}\beta _n^{(q)}(\lambda )u_{k-q}(X_n)$
and
\begin{equation}
\label{incrementMn}
u_k(M_{n+1})-u_k(M_n)
=\sum _{q=0}^{k-2}\frac 1{q!}\beta _{n+1}^{(q)}(\lambda )
\left[ u_{k-q}(X_{n+1})-
\frac{\beta _n^{(q)}(\lambda )}
{\beta _{n+1}^{(q)}(\lambda )}u_{k-q}(X_n)\right] .
\end{equation}
One can write
\begin{eqnarray*}
u_{k-q}(X_{n+1})-
\frac{\beta _n^{(q)}(\lambda )}{\beta _{n+1}^{(q)}(\lambda )}u_{k-q}(X_n)
&=&u_{k-q}(X_{n+1}-X_n)\\
&+&\left[ 1-\frac{\beta _n^{(q)}(\lambda )}
{\beta _{n+1}^{(q)}(\lambda )}\right] u_{k-q}(X_n).
\end{eqnarray*}
The relation
$\beta _n(\lambda )=(1+\lambda /(n+\tau _1-1))\beta _{n+1}(\lambda )$
implies, with Leibnitz formula, that
$$
1-\frac{\beta _n^{(q)}(\lambda )}{\beta _{n+1}^{(q)}(\lambda )}\in O(\frac 1n)
.$$
Besides, definition of the process $(X_n)_n$
(Definition \ref{defiPolya}) ensures that
$X_{n+1}-X_n\in\{w_1,\dots ,w_s\}$ is almost surely $O(1)$ and consequently
that $X_n$ is almost surely $O(n)$ as $n$ goes off to infinity
(elementary induction).
Hence
\begin{equation}
\label{truc}
u_{k-q}(X_{n+1})-
\frac{\beta _n^{(q)}(\lambda )}{\beta _{n+1}^{(q)}(\lambda )}u_{k-q}(X_n)
\in O(1)
\end{equation}
almost surely, as $n$ tends to infinity.
With the same tools as for the derivatives of $\gamma _{\tau _1,n}$
(see (\ref{derivGamma})), for every nonnegative integer $q$,
\begin{equation}
\label{beta}
\beta _n^{(q)}(\lambda )=(-1)^q
\frac{\log ^qn}{n^\lambda }
\frac{\Gamma (\tau _1+\lambda )}{\Gamma (\tau _1)}
+o\left( \frac{\log ^qn}{n^{\Re\lambda }}\right)
\end{equation}
as $n$ tends to infinity.
Thus (\ref{incrementMn}), (\ref{truc}) and (\ref{beta}) lead to
$$
u_k(M_{n+1})-u_k(M_n)
\in O\left(\frac{\log ^{k-2}n}{n^{\Re\lambda }}\right)
$$
almost surely as $n$ tends to infinity.
In particular, $|u_k(M_{n+1})-u_k(M_n)|^2$ is almost surely the general
term of a convergent series: the quadratic variation of the martingale
$(u_k(M_n))_n$ is almost surely bounded and the convergence part of the claim
is proved.

Almost surely and in ${\rm L}^p$ for every $p\geq 1$,
$$
\pi _JX_n=\gamma _{\tau _1,n}(\pi_JA)
\left[ \gamma _{\tau _1,n}(\pi _JA)^{-1}\pi _JX_n\right]
=\gamma _{\tau _1,n}(\pi_JA)(M_J+o(1))
$$
as $n$ tends to infinity.
As for equation (\ref{gamman-1}), one has
$$
\gamma _{\tau _1,n}(\pi _JA)
=\gamma _{\tau _1,n}(\lambda +N)\pi _J
=\sum _{k=2}^{\nu +2}
\left( \sum _{q=0}^{k-2}\frac 1{q!}
\gamma _{\tau _1,n}^{(q)}(\lambda )u_{k-q}\right) v_k
$$
and the asymptotics of the derivatives of $\gamma _{\tau _1,n}$
(see (\ref{derivGamma})) implies
$$
\pi _JX_n=
\frac {n^\lambda\log ^\nu n}{\nu !}
\frac {\Gamma (\tau _1)}{\Gamma (\tau _1+\lambda )}
u_2(M_J)v_{\nu +2}+o\left( n^\sigma\log ^\nu n\right) 
$$
which is the expected result (\ref{asymptpiJXn}) on $\pi _JX_n$.
Equation (\ref{gamman-1}) shows that
$u_2(M_n)=\beta _n(\lambda )u_2(X_n)$
and makes the proof of the claim complete with the help of (\ref{beta}).

\vskip 5pt
\noindent
$\bullet$
As in the proof of the large and principally semisimple case 
(Theorem \ref{lpss}),
$\pi _1X_n=(n+\tau _1-1)v_1$,
and the process splits into the sum
$$
X_n=nv_1+\sum _{k=2}^r\pi _{J_k}X_n
+Y_n+Z_n
$$
where
$Y_n=\sum _J\pi _JX_n$
the sum being extended to all monogenic blocks of indices different from any
$J_k$ that correspond to roots having real parts $>1/2$ and
$Z_n=\sum _{\{ k,~\Re\lambda _k\leq 1/2\} }\pi _kX_n$.
We study separately all terms of this decomposition.

\noindent
$\bullet$
Because of Theorem \ref{jointmoments} part {\it 1-},
as in the end of the proof of the large and principally semisimple case, 
$Z_n\in o(n^{\sigma _2}\log ^\nu n)$ almost surely and
in ${\rm L}^p$ for every $p\geq 1$ (remember that $\pi _kX_n=u_k(X_n)v_k$).

\noindent
$\bullet$
Every $J$ in the definition of $Y_n$ satisfies the assumption of the claim
with a root's real part $<\sigma _2$ or a cardinality $\leq \nu$.
Thus almost surely and in ${\rm L}^p$ for every $p\geq 1$,
$$
Y_n=o(n^{\sigma _2}\log ^\nu n)
$$
as $n$ tends to infinity.

\noindent
$\bullet$
For every $k\in\{ 2,\dots ,r\}$, $J_k$ satisfies the assumption of the claim
and if one denotes
$$
W_k=\lim _{n\to\infty}\frac{u_{\min J_k}(X_n)}{n^{\lambda (J_k)}},
$$
one obtains
$$
\pi _{J_k}X_n=
\frac 1{\nu !}n^{\lambda (J_k)}\log ^\nu nW_kv_{\max J_k}
+o\left( n^{\sigma _2}\log ^\nu n\right)
$$
almost surely and in ${\rm L}^p$ for every $p\geq 1$, which completes the proof
of (\ref{asymptlpnss}).
Note that $u_{\min J_k}$ is an eigenform of $A$ and that
$\gamma _{\tau _1,n}(\lambda (J_k))^{-1}u_{\min J_k}(X_n)$
is an ${\rm L}^{\geq 1}$-convergent complex-valued martingale.

\noindent
$\bullet$
Take any $\alpha _2,\dots ,\alpha _r\in \g Z_{\geq 0}$.
Then $\alpha =\sum _{2\leq k\leq r}\alpha _k\delta _{\min J_k}$
is a semisimple power of large projections, and
$$
E\left( \prod _{2\leq k\leq r}W_k^{\alpha _k}\right)
=\lim _{n\to\infty}\frac1{n^{\langle\alpha,\lambda\rangle}}E\uu ^\alpha (X_n)
=\frac{\Gamma (\tau _1)}{\Gamma (\tau _1+\langle\alpha ,\lambda\rangle )}
Q_\alpha (X_1)
$$
as can be deduced from Assertion~{\it 3-} in Theorem~\ref{jointmoments}.
This completes the proof of Theorem~\ref{lpnss}.
\QED

\masection{Remarks and examples}
\label{Exemples-sec}

\subsection{Some remarks}
\label{rem-subsec}

{\bf 1- Average case study of a P\'olya process}

If $(X_n)_n$ is a large P\'olya process, its asymptotic expectation can
readily be deduced from Theorems~\ref{lpss} and~\ref{lpnss}.
Without using the whole result, if $(X_n)_n$ is \emph{any} P\'olya process, one can
simply argue as follows.
Thanks to Relation~(\ref{projections}),
$
EX_n=\sum _{1\leq k\leq s}Eu_k(X_n).v_k
$.
If $J$ is any monogenic block of indices,
the subspace $\Vect\{u_k,~k\in J\}$ is $\Phi$-stable so that
Proposition~\ref{reducedmoments} (which is elementary) applies.
Hence $Eu_k(X_n)\in O(n^{\Re \lambda _k}\log ^{\nu _{\delta _k}}n)$
when $\lambda _k\neq 1$ and
$Eu_k(X_n)=nu _k(X_1)+O(1)$ when $\lambda _k=1$
(in order to directly apply Proposition~\ref{reducedmoments}, remember that
$u_k=Q_{\delta _k}$).
These facts imply the following result, already given in \cite{AthreyaKarlin}
when $1$ is simple root.

\begin{Prop}
\label{asymptEsperance}
If $\Pi _1=\sum _{\{ k,~\lambda _k=1\}}\pi _k$ denotes the projection on the
eigensubspace $\ker (A-1)$, then, as $n$ goes off to infinity,
$$
E(X_n)=n\Pi _1 (X_1)+O(n^\tau)
$$
where
$\tau =\max\left(
\{\Re (\lambda ),~\lambda\in\Sp (A),~\lambda\neq 1\}\cup\{0\}
\right)$.
\end{Prop}

\vskip 10pt\noindent
{\bf 2- Drift when $1$ is simple root}

When $1$ is a simple root of a P\'olya process $(X_n)_n$, the normalisation
$X_n/n$ converges almost surely and in ${\rm L}^{\geq 1}$ to the non random
vector $v_1$.
This can be deduced from Theorem~\ref{jointmoments} and
Decomposition~(\ref{trisection}), by arguments like in the end of
Theorem~\ref{lpss}'s proof.
This result is valid for small and large processes, without any
irreducibility-type condition (compare with \cite{JansonFunctional}).

\vskip 10pt\noindent
{\bf 3- Small P\'olya processes}

As it has been told in Section~\ref{intro}, a small \emph{irreducible}
P\'olya process has a Gaussian limit after normalisation (for a precise
meaning of the present notion of irreducibility and complete results,
see~\cite{AthreyaKarlin} and~\cite{JansonFunctional}).
When the irreducibility assumption is released, this normality fails down.
This fact can be explained by our treatment.
We illustrate it in details in dimension $2$.

Take the general two-dimensional P\'olya process and choose coordinates
such that the forms $l_k$ are the coordinates forms in $\g R^2$.
The matrix of the replacement endomorphism $A$ have then the form
$\left(
\begin{array}{cc}
1-a&b\\
a&1-b\\
\end{array}
\right)$
where $a$ and $b$ are nonnegative reals (with restrictive conditions
(\ref{tenable}) if at least one of them is $>1$).
The process is small whenever $a+b\geq 1/2$ because $\sigma _2=1-a-b$.
Let's assume for our example that $a+b>1/2$.
If one makes the choice $u_2=ax-by$, computation of the first reduced
polynomials shows that
\begin{equation}
\label{urn2general}
u_2^2=Q_{(0,2)}
-(a-b)(1-a-b)u_2
+\frac{ab(1-a-b)^2}{2(a+b)-1}u_1.
\end{equation}
The term of $Eu_2^2(X_n)$ having the highest order of magnitude is
$Eu_1(X_n)=nu_1(X_1)$, but its coefficient is zero if
$a$ or $b$ vanish or if $a+b=1$.
Such considerations justify the fact that the study of small triangular urns,
that are not irreducible,
has to be done separately in terms of asymptotics and limit laws
(see \cite{JansonFunctional}, \cite{JansonTriangular}, \cite{Puyhaubert}).

In arbitrary dimension $s$, one can refine the error term in
Assertion~{\it 1-} of Theorem~\ref{jointmoments}, but this
refinement requires more careful use of the replacement endomorphism.
This fact comes from the expansion
$E\uu ^\alpha(X_n)=EQ_\alpha (X_n)+
\sum _{\beta<\alpha }q_{\alpha ,\beta}EQ_\beta (X_n)$:
if $\alpha$ is a power of small projections, the term in the equality's second
member having
the highest order of magnitude as $n$ goes off to infinity is not necessarily
$EQ_\alpha (X_n)$, but $EQ_\alpha (X_n)$ \emph{may} nevertheless
be the winner if suitable coefficients $q_{\alpha ,\beta}$ vanish.

\vskip 10pt\noindent
{\bf 4- Limit random variables $W_k$}

As it can be seen in the proof of Theorem~\ref{lpss},
for any $k\in\{ 2,\dots ,r\}$,
the random variable $W_k$ is defined as the limit of the process
$u_k(X_n)/n^{\lambda _k}$ as $n$ tends to infinity.
This convergence is almost sure and in any ${\rm L}^p$, $p\geq 1$ and is
proved by martingale techniques.

To know whether $W_k$ is zero or not, it is sufficient to check the nullity of
$EW_k^2=\Gamma (\tau _1)Q_{2\delta _k}(X_1)/\Gamma (\tau _1+2\lambda _k)$
when
$W_k$ is real-valued (that is when $\lambda _k$ is real),
or of
$E|W_k|^2=\Gamma (\tau _1)Q_{\delta _k+\delta_{k'}}(X_1)/
\Gamma (\tau _1+2\Re\lambda _k)$
when $W_k$ is not real-valued ({\it i.e.} when
$\lambda _k\in\g C\setminus\g R$), where $k'$ is such that
$\overline{u_k}=u_{k'}$
(conditionally to the choice of a suitable Jordan basis).

Questions: what can be said about these variables?
Are the laws of the $W_k$ always determined by their moments?
Can they always be described in terms of known densities or other
distributions?

All these remarks and questions can readily be adapted to limit variables
$W_k$ of Theorem~\ref{lpnss}.

\vskip 10pt\noindent
{\bf 5- Conjugate replacement endomorphisms}

In the asymptotic almost sure expansions~(\ref{asymptlpss})
or~(\ref{asymptlpnss}), $\sigma _2$,
the complex numbers $\lambda _k$ and $\lambda (J_k)$ and the integer $\nu$
depend only on the conjugacy class of the replacement endomorphism $A$.
On the contrary, the distributions of the random variables $W_k$ depend on the
increment vectors $w_k$ and on the linear forms $l_k$
(and on initial condition $X_1$), but not only on the conjugacy class of
$A=\sum _kl_k\otimes w_k$: \emph{two processes having conjugate replacement
endomorphisms have the same asymptotic form (\ref{asymptlpnss}), but have
in general different limit laws $W_k$}.
For example, the two standardized large urns having
$R=\left(
\begin{array}{cc}
1&0\\
9/20&11/20\\
\end{array}
\right)$
and
$R'=\left(
\begin{array}{cc}
3/4&1/4\\
1/5&4/5\\
\end{array}
\right)$
as (conjugate) replacement matrices have respective second reduced polynomials
$Q_{(0,2)}=u_2(u_2+11/20)$ (see (\ref{triangle2})) and
$Q'_{(0,2)}={u'_2}^2-\frac{11}{400}u'_2+\frac{121}{800}u'_1$
(see (\ref{urn2general}), evident notations).
The algebraic relations satisfied by the  moments of $W$ and $W'$ are
not of the same kind.

Another way to formulate this remark, as suggested by the referee, is the
following.
Two processes may have the \emph{same} replacement endomorphism $A$
(which is the restriction of $\Phi$ over linear forms)
without having the same transition operator $\Phi$:
this will imply in general different $Q_\alpha$'s, even though the asymptotic
form will be of the same nature.
Note however that having the same replacement endomorphism $A$ does not mean
having the same linear forms $l_k$ and increment vectors $w_k$.

A natural question arises: when two processes have conjugate (or equal)
replacement endomorphisms, are their limit laws connected by some functional
relation?

\subsection{Examples}
\label{ex-subsec}

{\bf 1- P\'olya-Eggenberger urns}

\noindent
As stated in Section~\ref{intro}, any P\'olya-Eggenberger urn is
a P\'olya process after standardization, {\it i.e.} after division by $S$
in order to get balance equal to $1$.
For further developments of examples on the general two dimensional urn
process, on some generic examples in dimension $5$ and on the so-called
$s$-dimensional {\it cyclic urn} whose (semisimple) replacement matrix is
$$
\left(
\begin{array}{ccccc}
0&1&&&0\\
&0&1&&\\
&&0&&\\
&&&\ddots&1\\
1&&&&0
\end{array}
\right) ,
$$
see \cite{Barcelona}
(the cyclic urn defines a small P\'olya process if and only if $s\leq 6$
because $\sigma _2=\cos (2\pi /s)$).
In the present article, see (\ref{trigo2general}) for some developments on
the general triangular urn with two colours;
other considerations are made on the same subject in \cite{Puyhaubert}.

\vskip 10pt\noindent
{\bf 2- Triangular urns with two types of balls.}

\noindent
The general two-dimensional balanced triangular P\'olya urn
(generalized to real numbers) has the following $R$ as replacement matrix:
\begin{equation}
\label{trigo2general}
R=
\left(
\begin{array}{cc}
1&0\\
1-\ell &\ell\\
\end{array}
\right) ,
\end{equation}
where $\ell$ is any real number $\leq 1$.
In terms of P\'olya process, this means that $l_1$ and
$l_2$ are the coordinate forms, $w_1=\transp (1,0)$ and
$w_2=\transp (1-\ell ,\ell )$.
If one chooses $u_2(x,y)=y$ as second form for a Jordan basis,
a straightforward computation shows that
for any integer $p\geq 0$, one has
$Q_{p\delta _2}=u_2(u_2+\ell )\dots (u_2+(p-1)\ell )$
and $\Phi (Q_{p\delta _2})=p\ell Q_{p\delta _2}$
(the simple computation of the image by $\Phi$ of the product
$u_2(u_2+\ell )\dots$
suffices to show that this product equals $Q_{p\delta _2}$).
Reversing this last formula leads, for any integer $p\geq 0$, to
\begin{equation}
\label{triangle2}
u_2^p=\sum _{k=1}^p(-\ell )^{p-k}\stirling2 pkQ_{k\delta _2}
\end{equation}
where $\stirling2 pk$ denote Stirling numbers of the second kind
(see for example \cite{ConcreteMath} for this reversion formula).

In particular, if $\ell >0$, since the order of magnitude of
$EQ_{p\delta _2}(X_n)$ is $n^{p\ell}$ (Proposition \ref{reducedmoments}),
$Eu_2(X_n/n^\ell )^p$ tends to
$\ell ^p\times\Gamma (x_1+y_1)/\Gamma (x_1+y_1+p\ell )
\times\Gamma (y_1/\ell +p)/\Gamma (y_1/\ell )$
as $n$ tends to infinity, where $X_1=\transp (x_1,y_1)$
is the initial composition of the urn.
This shows the convergence in distribution of
$(X_n-nv_1)/n^\ell =u_2(X_n/n^\ell )v_2$
to the law having the written above expression as $p$-th moment
(the asymptotics of the computed $p$-th moment as $p$ tends to infinity
shows by means of Stirling formula that the limit law is determined by its
moments, proving the convergence in law; see for example \cite{Billingsley}
for relations between convergence of moments and convergence in distribution).
For descriptions of this limit laws in some very particular cases of parameters
$X_1$ and $\ell$ in terms of stable laws or Mittag-Leffler distribution, one
can refer to \cite{Puyhaubert} or \cite{JansonTriangular}.
When $\ell >1/2$, the process is large so that this convergence is almost sure
and in any ${\rm L}^p$, $p\geq 1$.

The case $\ell =0$ is degenerate: the process is deterministic.

When $\ell <0$, as $EQ_{p\delta _2}(X_n)\in O(n^{p\ell })$
(in any case, even if $\tau _1+p\ell$ is a nonpositive integer),
Formula~(\ref{triangle2}) implies that $Eu_2(X_n)^p=O(n^{\ell})$ for any $p$.
In this case, $u_2(X_n)$ tends almost surely to zero because balls of the
second type can never be added
(see for example~\cite{JansonTriangular}, Section 2, \emph{Degenerate cases}).

One can compare this to the results of \cite{Puyhaubert} and
\cite{JansonTriangular}.
It can easily be generalized to some classes of triangular urns of higher
dimension, principally semisimple or not (with enough zero entries, see
\cite{Barcelona} for examples).

\vskip 10pt\noindent
{\bf 3- Example of random replacement matrices}

The following example of urn process comes from a private communication of
Bernard Ycart.
Take an urn containing first $b$ black balls, $w$ white balls and one red ball.
As in the case of P\'olya urn processes, one draws successively balls from
the urn, with the following replacement rule.
If a black (respectively white) ball is drawn, replace it in the urn together
with another black (resp. white) one.
If the red ball is drawn, replace it in the urn together with a black one with
probability $p\in [0,1]$ or a white one with probability $1-p$.

As it is described, this urn process is not P\'olya.
But it is equivalent to the P\'olya process defined in $\g R^4$ by:
the $l_k$ are the coordinate forms, the replacement matrix
({\it i.e.} the matrix whose rows are the coordinates of the $w_k$) is
$$R=
\left(
\begin{array}{cccc}
1&0&0&0\\
0&1&0&0\\
1&0&0&0\\
0&1&0&0\\
\end{array}
\right),
$$
and the initial vector is $(b,w,p,1-p)$.
It can be viewed as ``non-integer'' four-colour P\'olya-Eggenberger urn
process, the colours being black, white, dark red and light red, the
replacement matrix being $R$.
Only the non-integer initial vector prevents our first problem from being 
a true P\'olya-Eggenberger urn process.
The matrix $R$ admits $1$ as double root, so that $X_n/n$ converges almost
surely and its limit has Dirichlet distribution (see example {\bf 7-}).

This example can easily be generalised to other replacement rules, provided
that one never adds any red ball.

\vskip 10pt\noindent
{\bf 4- $m$-ary search trees}

``$m$-ary search trees are fundamental data structures in computer science
used in searching and sorting '' (citation from~\cite{FillKapur}).
The space-requirements vector of an $m$-ary search tree under the random
permutation model is an $m-1$-dimensional P\'olya process as can be seen
in~\cite{ChauvinPouyanne}.
It only appears under the form of an urn process after some suitable change
of coordinates.
The associated endomorphism $A$ is semisimple and the process is large if,
and only if $m\geq 27$.
One can find further developments on this large process in \cite{Barcelona}.
See \cite{ChernHwang}, \cite{JansonFunctional} and \cite{FillKapur}
for different treatments of the subject.

\vskip 10pt\noindent
{\bf 5- Random $2-3$-trees}

This example comes from data structures in computer science too.
The repartition of external nodes of a random $2-3$-tree having $1$ or $2$
sisters is the two-dimensional P\'olya-Eggenberger urn process with
initial condition $X_1=\transp (2,0)$ and replacement matrix
$\left(
\begin{array}{cc}
-2&3\\
4&-3\\
\end{array}
\right)$.
This process is small ($\sigma _2=-6$) and principally semisimple.
It follows from \cite{JansonFunctional} that its second order term has
normal distribution.
This example is the base example of \cite{FGP}.

If one goes one step further, one can distinguish external nodes of a
random $2-3$ tree with regard to the shape of the descendants-tree
of their grand-mothers.
This process is a $10$-dimensional urn process with balance $1$.
Its replacement matrix
\def\st{\scriptstyle}
$$R=\left(
\begin{array}{cccccccccc}
\st -4&\st 2&\st 3&\st 0&\st 0&\st 0&\st 0&\st 0&\st 0&\st 0\\
\st 0&\st -2&\st -3&\st 6&\st 0&\st 0&\st 0&\st 0&\st 0&\st 0\\
\st 0&\st -2&\st -3&\st 0&\st 6&\st 0&\st 0&\st 0&\st 0&\st 0\\
\st 0&\st 0&\st 0&\st -6&\st 0&\st 4&\st 3&\st 0&\st 0&\st 0\\
\st 0&\st 0&\st 0&\st 0&\st -6&\st 4&\st 3&\st 0&\st 0&\st 0\\
\st 0&\st 0&\st 0&\st \st 0&\st 0&\st -4&\st -3&\st 2&\st 6&\st 0\\
\st 8&\st 0&\st 0&\st 0&\st 0&\st -4&\st -3&\st 0&\st 0&\st 0\\
\st 0&\st 0&\st 0&\st 0&\st 0&\st 0&\st 0&\st -2&\st -6&\st 9\\
\st 4&\st 2&\st 3&\st 0&\st 0&\st 0&\st 0&\st -2&\st -6&\st 0\\
\st 4&\st 0&\st 0&\st 6&\st 0&\st 0&\st 0&\st 0&\st 0&\st -9\\
\end{array}
\right)
$$
contains negative off-diagonal entries.
This does not prevent the urn to be tenable
(for ``physical'' reasons, as first argument!).
Indeed, the columns of $R$ containing
these negative entries are coupled in the following sense:
if $j\neq k$ and $r_{j,k}<0$, then the columns of $r_{j,k}$ and $r_{k,j}$
are proportional.
These proportionalities imply deterministic relations between number of balls
of concerned colours.
For example, at any time, the number of (algebraically) added balls of colour
$9$ is thrice the number of added balls of colour $8$ so that when a ball
of colour $8$ is drawn, if one can subtract $2$ balls of colour $8$,
one can subtract $6$ balls of colour $9$ as well.
The same kind of property holds for balls of colours $2$ and $3$, and
for balls of colours $6$ and $7$.
For such reasons, the same recurrence that shows that a P\'olya process does
not extinguish shows that our urn is tenable.

Moreover, our treatment of P\'olya processes readily applies to this urn
process.
It is small and principally semisimple, with $\sigma _2=0$ (the multiplicity
of the eigenvalue $0$ of $R$ is $3$).
Its study shows for instance that, if $n$ is the number of external nodes
of the tree, the average number of their grand-mothers is $\sim 0.182n$,
that on average $\sim 21\%$ (resp. $\sim 24\%$) of external nodes have
grand-mothers having themselves $4$ (resp. $5$) grand-children {\it etc.}

Patient readers can go still one step further, looking at the
fourth level of genealogical trees of external nodes.
This leads to the study of a $76$-dimensional urn process.

\vskip 10pt\noindent
{\bf 6- Congruence in binary search trees}
  
\noindent
The following example is mentioned in \cite{ChernFuchsHwang} as a private
unpublished idea of S. Janson\footnote{
S. Janson has developed his example in~\cite{JansonCongruences}
during the revision of the present article.
}.
Take a binary search tree and an integer $s\geq 2$.
Consider the random vector of $\g R^s$ whose $k$-th coordinate is the number
of leaves whose depth is $\equiv k~[mod~s]$.
This defines an $s$-colour urn process with (semisimple) replacement matrix
$$
\left(
\begin{array}{ccccc}
-1&2&&&\\
&-1&2&&\\
&&-1&&\\
&&&\ddots&2\\
2&&&&-1
\end{array}
\right) ;
$$
the balance is one and $\sigma _2=-1+2\cos (2\pi /s)$, so that the urn is
small if and only if $s\leq 8$.
As it is readily irreducible, it can be deduced from \cite{JansonFunctional}
that its second order term has normal distribution when $s\leq 8$.
When $s\geq 9$, the process is large and its asymptotics is described by
Theorem \ref{lpss}.

\vskip 10pt\noindent
{\bf 7- Processes having $1$ as multiple root}

Let $(X_n)_n$ be a P\'olya process having $1$ as multiple root;
the way to use Theorem~\ref{lpss}
to determine the almost sure limit law of $X_n/n$ suggests
to abandon our convention $u_1=\sum _{k=1}^sl_k$.
This does not change the validity of the whole result.

Let $r\geq 2$ be the multiplicity of $1$ as eigenvalue of $A$.
We choose a basis $(u_1,\dots ,u_r)$ of $A$-fixed linear forms ({\it i.e.} a
basis of $\ker (^t\! A-1)$), using the classical following construction.
Consider the graph $\cal G$ whose vertices are the numbers $\{ 1,\dots ,s\}$
and where two vertices $i$ and $j$ are connected by an edge when
$l_i(w_j)\neq 0$ or $l_j(w_i)\neq 0$.
Let $I_1,\dots ,I_r$ be the connected components of $\cal G$ (the fact that
there are $r$ such components in a consequence of what follows).
For any $(j,k)$, it is readily shown that $l_j(w_k)=0$ if
$j\notin I_k$.
We define
$$u_k=\sum _{j\in I_k}l_j$$
for any $k\in\{ 1,\dots ,r\}$, so that $u_k(w_j)=0$ if $j\notin I_k$
and $u_k(w_j)=\sum _{1\leq i\leq s}l_i(w_j)=1$ if $j\in I_k$.
A straightforward computation shows that any $u_k$ is an $A$-fixed linear form.
Moreover, the restriction of $^t\! A$ to the stable subspace spanned by the
$l_j$, $j\in I_k$ is irreducible so that, because of Perron-Frobenius theory,
$u_k$ spans the unique line of $A$-fixed forms of this subspace.
This shows that $(u_1,\dots ,u_r)$ is a basis of $A$-fixed linear forms.
This basis is then completed into a Jordan basis $(u_1,\dots ,u_s)$ under
conditions {\it 2-} of Definition~\ref{jordanBasis-def}.

For such a basis, $\sum _{k=1}^sl_k=\sum _{k=1}^ru_k$.
The properties of $u_k$'s imply in particular that
for any $\alpha =(\alpha _1,\dots ,\alpha _r,0,\dots )$,
$$
Q_\alpha =\prod _{k=1}^ru_k(u_k+1)\dots (u_k+\alpha _k-1)
$$
and that $Q_\alpha$ is eigenfunction for $\Phi$, associated with the eigenvalue
$|\alpha |=\sum _{k=1}^r\alpha _k$.
It follows then from Theorem \ref{lpss} that $X_n/n$ converges almost
surely and in any ${\rm L}^p$, $p\geq 1$ to a random vector
$\sum _{k=1}^rW_kv_k$, where the joint moments of the real random variables
$W_1,\dots ,W_r$ are given by
$$
EW^\alpha =\frac{\Gamma (\tau _1)}{\Gamma (\tau _1+|\alpha |)}
\prod _{k=1}^r\frac{\Gamma (u_k(X_1)+\alpha _k)}{\Gamma (u_k(X_1))}.
$$
One recognizes here the moments of a Dirichlet distribution with parameters
$u_1(X_1),\dots ,u_r(X_1)$ whose density on the simplex
$\{ x_1\geq 0,\dots ,x_r\geq 0,\sum _{k=1}^rx_k=1\}$
of $\g R^r$ is given by
$$
(x_1,\dots ,x_r)\mapsto
\Gamma \left(\sum _{k=1}^ru_k(X_1)\right)
\prod _{k=1}^r\frac{x_k^{u_k(X_1)}}{\Gamma (u_k(X_1))}
$$
(see \cite{Gouet}).
This distribution is obviously characterized by its moments.
In reference to the original paper of P\'olya, processes under this assumption
have been called {\bf essentially P\'olya} in~\cite{Barcelona}.

\begin{small}

\end{small}

\end{document}